\documentclass{amsart}
\usepackage{mathrsfs}
\usepackage{amssymb}
\usepackage{amsmath,amssymb}

\newtheorem{theorem}{Theorem}[section]
\newtheorem{lemma}[theorem]{Lemma}

\theoremstyle{definition}

\newtheorem{corollary}[theorem]{Corollary}
\theoremstyle{definition}

\theoremstyle{remark}
\newtheorem{remark}[theorem]{Remark}

\numberwithin{equation}{section}

\begin{document}

\title[Carleson measure problems
for  parabolic Bergman spaces]{Carleson measure problems for
parabolic Bergman spaces  and homogeneous  Sobolev spaces}



\author{Zhichun Zhai}
\address{Department of Mathematics and Statistics, Memorial University of Newfoundland, St. John's, NL A1C 5S7, Canada}
\curraddr{}
 \email{a64zz@mun.ca}
\thanks{Project supported in part  by Natural Science and
Engineering Research Council of Canada.}



\date{}

\dedicatory{}

\begin{abstract}
 Let  $b_{\alpha}^{p}(\mathbb{R}^{1+n}_{+})$  be the space of  solutions
 to the  parabolic equation $\partial_{t}u+(-\triangle)^{\alpha}u=0$
 $(\alpha\in(0, 1])$ having finite $L^{p}(\mathbb{R}^{1+n}_{+})$ norm.
 We  characterize  nonnegative Radon measures $\mu$ on $\mathbb{R}^{1+n}_{+}$
  having the property $\|u\|_{L^{q}(\mathbb{R}^{1+n}_{+},\mu)}\lesssim
\|u\|_{\dot{W}^{1,p}(\mathbb{R}^{1+n}_{+})},$ $1\leq p\leq
q<\infty,$ whenever $u(t,x)\in
b_{\alpha}^{p}(\mathbb{R}^{1+n}_{+})\cap
\dot{W}^{1.p}(\mathbb{R}^{1+n}_{+}).$
 Meanwhile,  denoting by $v(t,x)$  the solution of the above equation with Cauchy data
  $v_{0}(x),$ we characterize nonnegative Radon  measures $\mu$ on $\mathbb{R}_{+}^{1+n}$ satisfying
$\|v(t^{2\alpha},x)\|_{L^{q}(\mathbb{R}_{+}^{1+n},
\mu)}\lesssim\|v_{0}\|_{\dot{W}^{\beta,p}(\mathbb{R}^{n})},$
$\beta\in (0,n),$ $p\in [1, n/\beta],$ $q\in(0, \infty).$  Moreover,
we obtain the
decay of $v(t,x),$  an iso$-$capacitary inequality and  a trace inequality.\vspace{0.1in}\\
\end{abstract}

\maketitle

\tableofcontents \pagenumbering{arabic}

 {\section{Introduction and Statement of  Results}}

 Carleson measure was first introduced   in classical  Hardy space
 (see Carleson \cite{L. Carleson}) and have been extensively studied,
 for example, see Dafni-Karadzhov-Xiao \cite{Dafni karadzhov Xiao},
 Dafni-Xiao \cite{Dafni Xiao 2},  Hastings \cite{W. Hastings},
 Johnson \cite{R. Johnson}, and Xiao \cite{J. Xiao}-\cite{J. Xiao 3}
 and the  references therein. This paper   considers  Carleson measure
   problems via the parabolic equation
\begin{equation}\label{a7}
 \partial_{t}u(t,x)+(-\triangle)^{\alpha}u(t,x)=0, (t,x)\in\mathbb{R}^{1+n}_{+}
 =(0,\infty)\times \mathbb{R}^{n}
\end{equation}
 with $\alpha\in(0,1],$  and its Cauchy problem
\begin{equation}\label{a1}
 \left
 \{\begin{array}{ll}
 \partial_{t}v+(-\triangle)^{\alpha}v=0, & (t,x)\in \mathbb{R}_{+}^{1+n}; \\
 v(0,x)=v_{0}(x), & x\in \mathbb{R}^{n}, \end{array} \right.
\end{equation}
    where $\bigtriangleup$ is the Laplacian with respect to $x$ and
$$
 (-\triangle)^{\alpha}u(t,x)
  =\mathcal {F}^{-1}(|\xi|^{2\alpha}\mathcal{F}(u(t,\xi)))(x)
 $$
  with  $\mathcal {F}$ and $\mathcal {F}^{-1}$  being  the Fourier  transform
  and the inverse Fourier transform. More specifically, we  characterize
  nonnegative Radon  measures $\mu$ on $\mathbb{R}^{1+n}_{+}$ having the property
 \begin{equation}\label{em1}
  \|u(t,x)\|_{L^{q}(\mathbb{R}^{1+n}_{+},\mu)}\lesssim\|\nabla_{(t,x)} u(t,x)\|_{L^{p}(\mathbb{R}^{1+n}_{+})},\  \forall
  u\in b_{\alpha}^{p}(\mathbb{R}^{1+n}_{+})\cap \dot{W}^{1,p}(\mathbb{R}^{1+n}_{+}),
\end{equation}
for  $1\leq p\leq q<\infty,$ or
 \begin{equation}\label{em2}
  \|v(t^{2\alpha},x)\|_{L^{q}(\mathbb{R}_{+}^{1+n}, \mu)}\lesssim
  \|v_{0}(x)\|_{\dot{W}^{\beta,p}(\mathbb{R}^{n})}, \forall v_{0}\in
  \dot{W}^{\beta,p}(\mathbb{R}^{n}),
\end{equation}
 for  $\beta\in (0,n),$ $p\in [1, n/\beta]$ and  $q\in(0, \infty).$ Here
 $$
 v(t,x)=S_{\alpha}(t)v_{0}(x):=K^{\alpha}_{t}(x)\ast v_{0}(x)
 $$
solves  (\ref{a1}),
$$K_{t}^{\alpha}(x)=(2\pi)^{-\frac{n}{2}}\int_{\mathbb{R}^{n}}e^{i
x\cdot \xi}e^{-t|\xi|^{2\alpha}}d\xi\geq 0,\ \forall (t,x)\in
\mathbb{R}^{1+n}_{+}
$$
and $g(x)\ast h(x)$ means  the convolution between $g(x)$ and $f(x)$
on the space variable.

The  main motivation of considering embeddings (\ref{em1}) and
(\ref{em2}) comes from  the so called trace inequalities problem.
  Particularly, for a nonnegative Borel measure $\mu$ on $\mathbb{R}^{n},$ when working
   on spectral problems for Schr¡§odinger operators, Maz'ya
first discovered in 1962 (see Maz'ya \cite{V.G. Maz'ya
62}-\cite{V.G. Maz'ya 64} and \cite{V.G. Maz'ya 1}) that  if
$1<p\leq q$ and  $pl<n$ then
\begin{equation}
\label{equavi1} \|u\|_{L^{q}(\mathbb{R}^{n},\mu)}\lesssim
\|u\|_{h_{p}^{l}(\mathbb{R}^{n})}, \forall u\in
h_{p}^{l}(\mathbb{R}^{n})
\end{equation}
holds if and only if
\begin{equation} \label{equavi2}
\sup\left\{\frac{(\mu(E))^{p/q}}{cap(E, h_{p}^{l})}: E\subset
\mathbb{R}^{n}, cap(E, h_{p}^{l})>0\right\}<\infty.
\end{equation}
 Here  $U\lesssim V$ denotes $U\leq \theta V$
for some  positive $\theta$ which is independent  of the sets or
functions under consideration in both
 $U$ and $V,$
 $h_{p}^{l}(\mathbb{R}^{n})$
is the completion of $C_{0}^{\infty}(\mathbb{R}^{n})$ with respect
to
 $$
 \|f\|_{h_{p}^{l}(\mathbb{R}^{n})}=\|(-\triangle)^{l/2}f\|_{L^{p}}
  $$
  and  $cap(E, h_{p}^{l}(\mathbb{R}^{n}))$ is the capacity of
$E$ associated with $h_{p}^{l}(\mathbb{R}^{n}).$  Such embeddings
like (\ref{equavi1}) are referred to as trace inequalities, see
Adams-Hedberg \cite{D.R. Adams Hedberg}. Meanwhile, (\ref{equavi2})
is  called  isocapacitary inequality, see Maz'ya \cite{V.G. Maz'ya
4}. Since Maz'ya established the pioneer work in \cite{V.G. Maz'ya
62}-\cite{V.G. Maz'ya 64}, other equivalent conditions of trace
inequalities 
were established
 by Maz'ya \cite{V.G. Maz'ya 2}-\cite{Maz'ya},
  Maz'ya-Preobra$\breve{z}$enski$\breve{{\i}}$ \cite{Maz'ya
  Preobra}, Maz'ya-Verbitsky \cite{Maz'ya Verbitsky 1}, Adams \cite{D.R. Andams 76}
 and new advances of such problems were made by
 Cascante-Ortega-Verbitsky \cite{C. CASCANTE J.M. ORTEGA AND I.E.
VERBITSKY 1} in which they established similar trace inequality for
a general class of radially decreasing convolution kernels. When
$0<q<p$ and $1<p<\infty,$ (\ref{equavi1}) holds if and
  only if
\begin{equation}\label{equavi4}
 \int_{0}^{\infty}\left(\frac{t^{p}}{\vartheta(t)^{q}}\right)^{1/(p-q)}\frac{dt}{t}<\infty
 \end{equation}
 with $\vartheta(t)=\inf\{cap(E, S_{p}^{l}): {E\subset \mathbb{R}^{n}, \mu (E)\geq
 t}\},$ see Maz'ya-Netrusov \cite{Mazya Netrusov} and Verbitsky \cite{Verbitsky}.
When $1<p<q<\infty,$ (\ref{equavi1}) holds if and only if
\begin{equation}\label{eqnavi3}
\sup_{x\in \mathbb{R}^{n},
r>0}\frac{(\mu(B(x,r)))^{p/q}}{r^{(n-lp)}}<\infty,
\end{equation}
 where $B(x,r)$
is a ball of radius $r$ centered at $x\in \mathbb{R}^{n},$ see
Adams-Hedberg \cite[Theorem 7.2.2]{D.R. Adams Hedberg}.  When
$0<q<p$ and  $1<p<n/l,$ (\ref{equavi1}) holds if and only if the
Wolff potential
$$
W_{\alpha,p}^{\mu}(x):=\int_{0}^{\infty}(r^{lp-n}\mu(B(x,r)))^{p'-1}\frac{dt}{t}\in L^{q(p-1)/(p-q)}(\mu),
$$
see, Cascan-Ortega-Verbitsky \cite{C. CASCANTE J.M. ORTEGA AND I.E.
VERBITSKY 1}.  There exist  other
 conditions involving no capacity, which are equivalent to (\ref{equavi1}), see, for example,
  Maz'ya \cite{Maz'ya 5}, Maz'ya-Verbitsky \cite{Maz'ya Verbitsky 1}
 and Verbitsky \cite{Verbitsky}.   These equivalent conditions were
 widely applied to harmonic analysis, operator theory, function
  spaces, linear and
 nonlinear partial differential equations, etc., see, Adams-Hedberg \cite{D.R. Adams
 Hedberg},
 Maz'ya \cite{V.G. Maz'ya 1} and Maz'ya-Shaposhnikova \cite{Mazya shaposhnikova} and the references therein.

 This paper  characterizes (\ref{em1}) or (\ref{em2}) by conditions like
(\ref{equavi2}), (\ref{equavi4}) and (\ref{eqnavi3}). To do this, we
need  the following preliminary  materials.

     We   always assume
 that $\beta\in (0,n)\backslash\mathbb{N}$ when  $p=1$ or $n/\beta.$
$\dot{W}^{1,p}(\mathbb{R}^{1+n}_{+})$ is the completion  of
  $C^{\infty}_{0}(\mathbb{R}^{1+n}_{+})$ 
  with respect to the norm
 $$
 \|f\|_{\dot{W}^{1,p}(\mathbb{R}^{1+n}_{+})}=\left(\int_{\mathbb{R}^{1+n}_{+}}|\nabla_{(t,x)}f|^{p}dtdx\right)^{1/p}.
 $$
$b^{p}_{\alpha}(\mathbb{R}^{1+n}_{+}) (\alpha\in (0, 1])$ introduced
by  Nishio-Shimomura-Suzuki \cite{M. Nishio K. Shimomura N. Suzuki}
is the parabolic Bergamn space on $\mathbb{R}_{+}^{1+n},$ which is
the set of all solutions  of the  parabolic equation (\ref{a7})
having
 finite ${L^{p}(\mathbb{R}^{1+n}_{+})} $ norm.
 $\dot{W}^{\beta, p}(\mathbb{R}^{n})$
 is the homogeneous  Sobolev space which is the completion of 
 $C^{\infty}_{0}(\mathbb{R}^{n})$ with respect to the
norm
  \begin{equation}
\|f\|_{\dot{W}^{\beta,p}(\mathbb{R}^{n})} =\left
\{\begin{array}{lll} \|(-\triangle)^{\beta/2}f\|_{L^{p}}, & p\in
(1,n/\beta),\vspace{0.1in} \\
\left(\int\limits_{\mathbb{R}^{n}}\frac{\|\triangle^{k}_{h}f\|_{L^{p}}^{p}}{|h|^{n+p\beta}}dh\right)^{1/p},
& p=1\  \hbox{or}\  p=n/\beta\ , \beta\in (0,n)\backslash\mathbb{N},
\end{array} \right.\nonumber
  \end{equation}
 where
  \begin{equation}
\triangle^{k}_{h}f(x)=\left \{\begin{array}{lll}
\triangle^{1}_{h}\triangle^{k-1}_{h}f(x), & k>1, \vspace{0.1in}\\
f(x+h)-f(x), & k=1,
\end{array} \right.\nonumber
  \end{equation}
$k=1+[\beta],$ $\beta=[\beta]+\{\beta\}$ with $\{\beta\}\in (0,1).$

 If $X=\mathbb{R}^{1+n}_{+},$ $\beta=1$  and $p\geq 1,$ or $X=\mathbb{R}^{n},$ $\beta\in(0,n)$ and $p\in[1, n/\beta],$
 $cap_{\dot{W}^{\beta,p}(X)}(S)$ (see Maz'ya   \cite{V.G. Maz'ya 1}) is  the variational capacity  of an arbitrary set $S\subseteq X:$
 $$
 cap_{\dot{W}^{\beta,p}(X)}(S)=\inf\left\{\|f\|^{p}_{\dot{W}^{\beta,p}(X)}:
 f\in V_{X}(S) \right\}.
 $$
  Here
  $$
  V_{\mathbb{R}^{1+n}_{+}}(S)=\{f\in\dot{W}^{1,p}(\mathbb{R}^{1+n}_{+}):  S\subseteq
 \hbox{Int}(\{x\in \mathbb{R}^{1+n}_{+}: f\geq 1\})\}
 $$
and
 $$
  V_{\mathbb{R}^{n}}(S)=\{f\in\dot{W}^{\beta,p}(\mathbb{R}^{n}): f\geq 0, S\subseteq
 \hbox{Int}(\{x\in \mathbb{R}^{n}: f\geq 1\})\}
 $$
 with $\hbox{Int}(E)$ be  the interior of a set $E\subseteq X.$  For $t\in (0,\infty),$ $c_{p}^{\beta}(\mu;t)$  is the
 $(p,\beta)-$variational capacity minimizing function   associated with both  $\dot{W}^{\beta,p}(\mathbb{R}^{n})$ and
 a nonnegative measure $\mu$ on  $\mathbb{R}_{+}^{1+n}$  defined by
 $$
 c_{p}^{\beta}(\mu;t)=\inf\{cap_{\dot{W}^{\beta,p}(\mathbb{R}^{n})}(O):\ \hbox{bouded open}\ O\subseteq \mathbb{R}^{n},
 \mu(T(O))>t\},
 $$
where $T(O)$ is the tent based on an open subset $O$ of
$\mathbb{R}^{n}:$
 $$
 T(O)=\{(r,x)\in \mathbb{R}^{1+n}_{
 +}: B(x,r)\subseteq O\},
 $$
 with $B(x,r)$ be the open ball centered at $x\in \mathbb{R}^{n}$ with radius $r>0.$

For handling  the endpoint case $p=n/\beta$ we also need the
definition of the Riesz potentials (see  Adams-Xiao \cite{D. R.
Adams J. Xiao} and Adams \cite{D.R. Adams 4}) on $\mathbb{R}^{2n}$
as follows. The Riesz potential of order $\gamma\in(0,2n)$ is
defined by
 $$
 I_{\gamma}^{(2n)}\ast f(z)=\int_{\mathbb{R}^{2n}}|z-y|^{\gamma-2n}f(y)dy, \ z\in \mathbb{R}^{2n}.
 $$
 From Adams \cite[ Theorem 5.1]{D.R. Adams 4}, we have that if $u(x)$ and $I_{\gamma}^{(2n)}\ast |f|(x,0)$ are both in
 $L^{1}_{loc}(\mathbb{R}^{n})$  with
  \begin{equation}\label{a2}
   f(x,h)=|h|^{-\gamma}\triangle^{k}_{h}u(x),
   \end{equation}
 then  $u(x)=CI_{\gamma}^{(2n)}\ast f(x,0),$ for  a.e.  $x\in \mathbb{R}^{n}$ and some $C>0.$ Note that if
 $u\in \dot{W}^{\beta,n/\beta}(\mathbb{R}^{n})$ and $\gamma=2\beta\in(0,2n)$ then the function $f(\cdot,\cdot)$ in
(\ref{a2}) belongs to the space $L^{n/\beta}(\mathbb{R}^{2n}).$
For any $\gamma\in (0,2n),$
$\dot{\mathcal{L}}_{\gamma}^{p}(\mathbb{R}^{2n})=I_{\gamma}^{2n}\ast
L^{p}(\mathbb{R}^{2n})$ defined by   $\|I_{\gamma}^{2n}\ast
f\|_{\dot{\mathcal{L}}_{\gamma}^{p}(\mathbb{R}^{2n})}=\|f\|_{L^{p}_{2n}}.$

To state our main results, let us agree to  more  conventions.
$U\approx V$ if $U\lesssim V$ and $V\lesssim U;$ for $0< p,
q<\infty$ and a nonnegative  Radon measure $\mu$ on
$X=\mathbb{R}^{1+n}_{+}\ \hbox{or}\  \mathbb{R}^{n},$ $L^{q,p}(X,
\mu)$ and $L^{q}(X, \mu)$ denote the Lorentz space and the Lebesgue
space of all functions $f$ on $X$ for which
$$
\|f\|_{L^{q,p}(X, \mu)}=\left(\int_{0}^{\infty}(\mu(\{x\in X:
|f(x)|>\lambda\}))^{p/q}d\lambda^{p}\right)^{1/p}<\infty
$$
 and
$$
\|f\|_{L^{q}(X,
\mu)}=\left(\int_{X}|f(x)|^{q}d\mu\right)^{1/q}<\infty,
$$
 respectively. Moreover, we use $L^{q,\infty}(X, \mu)$ as the set of all $\mu-$measurable functions $f$ on $X$ with
$$
\|f\|_{L^{q,\infty}(X, \mu)}=\sup_{\lambda>0}\lambda(\mu(\{x\in X:
|f(x)|>\lambda\}))^{1/q}<\infty.
$$

\begin{theorem}\label{th1} Let   $1\leq p\leq q<\infty$  and $\mu$ be a nonnegative Radon  measure
 on $\mathbb{R}^{1+n}_{+}.$  Then the following statements are equivalent:\\
(a)
$$
\|u\|_{L^{q,p}(\mathbb{R}^{1+n}_{+}, \mu)}\lesssim\|u\|_{\dot{W}^{1,
p}(\mathbb{R}^{1+n}_{+})}, \ \ u\in \dot{W}^{1,
p}(\mathbb{R}^{1+n}_{+})\cap b^{p}_{\alpha}(\mathbb{R}^{1+n}_{+}),
$$
 (b)
$$
\|u\|_{L^{q}(\mathbb{R}^{1+n}_{+},
\mu)}\lesssim\|u\|_{\dot{W}^{1,p}(\mathbb{R}^{1+n}_{+})}, \ \ u\in
\dot{W}^{1, p}(\mathbb{R}^{1+n}_{+})\cap
b^{p}_{\alpha}(\mathbb{R}^{1+n}_{+}),
$$
 (c)
 $$
 \|u\|_{L^{q,\infty}(\mathbb{R}^{1+n}_{+}, \mu)}\lesssim\|u\|_{\dot{W}^{1, p}(\mathbb{R}^{1+n}_{+})}, \ \
u\in \dot{W}^{1, p}(\mathbb{R}^{1+n}_{+})\cap
b^{p}_{\alpha}(\mathbb{R}^{1+n}_{+}),
$$
 (d)
 $$
 (\mu(O))^{p/q}\lesssim cap_{\dot{W}^{1, p}(\mathbb{R}^{1+n}_{+})}(O), \
 \hbox{ open}\ O\subseteq \mathbb\mathbb{R}^{1+n}_{+}.
 $$
 If $0<q<p=1,$ then  $(b)\Longrightarrow(c)\Longrightarrow(d)\Longrightarrow(a).$
  \end{theorem}

 In the following,  $v(t,x)$ is the solution of equation (\ref{a1}) with Cauchy data $v_{0}(x).$

\begin{theorem} \label{th2} Let  $\beta\in(0,n),$ $0<q<p,$ $1< p\leq n/\beta$ and $\mu$  a nonnegative Radon measure on
$\mathbb{R}_{+}^{1+n}.$ Then the following two conditions are equivalent:\\
(a)
$$
\|v(t^{2\alpha},x)\|_{L^{q}(\mathbb{R}_{+}^{1+n}, \mu)}\lesssim
\|v_{0} \|_{\dot{W}^{\beta,p}(\mathbb{R}^{n})},\ \forall v_{0}\in
\dot{W}^{\beta, p}(\mathbb{R}^{n}).
$$
 (b)
 $$
 \int_{0}^{\infty}\left(\frac{t^{p/q}}{c_{p}^{\beta}(\mu; t)}\right)^{q/(p-q)}\frac{d t}{t}<\infty.
 $$
\end{theorem}

If we change   $1<p\leq n/\beta$ and $0<q<p$ into $1\leq p\leq
n/\beta$ and $p\leq q<\infty,$ then the conditions (a) and (b) of
Theorem \ref{th2} can be replaced by a weak-type one and two simpler
ones, respectively.

\begin{theorem} \label{th3}Let  $\beta\in(0,n),$   $1\leq p\leq n/\beta,$ $p\leq q<\infty$ and
$\mu$ a nonnegative Radon measure on $\mathbb{R}^{1+n}_{+}.$ Then the following five conditions are equivalent:\\
(a)
$$
\|v(t^{2\alpha},x)\|_{L^{q,p}(\mathbb{R}^{1+n}_{+}, \mu)}\lesssim
\|v_{0}\|_{\dot{W}^{\beta,p}(\mathbb{R}^{n})}, \ \forall v_{0}\in
\dot{W}^{\beta, p}(\mathbb{R}^{n}).
$$
(b)
$$
\|v(t^{2\alpha},x)\|_{L^{q}(\mathbb{R}_{+}^{1+n}, \mu)}\lesssim
\|v_{0}\|_{\dot{W}^{\beta,p}(\mathbb{R}^{n})}, \ \forall v_{0}\in
\dot{W}^{\beta, p}(\mathbb{R}^{n}).
$$
 (c)
 $$
 \|v(t^{2\alpha},x)\|_{L^{q,\infty}(\mathbb{R}^{1+n}_{+}, \mu)}\lesssim\|v_{0}\|_{\dot{W}^{\beta,p}(\mathbb{R}^{n})}, \
 \forall v_{0}\in \dot{W}^{\beta, p}(\mathbb{R}^{n}).
 $$
 (d)
 $$
 \sup\limits_{t>0}\frac{t^{p/q}}{c_{p}^{\beta}(\mu; t)}<\infty.
 $$
 (e)
 $$
 \sup\left\{\frac{\left(\mu\left(T(O)\right)\right)^{p/q}}
 {cap_{\dot{W}^{\beta,p}(\mathbb{R}^{n})}(O)}:\ \hbox{bounded open }\ O\subseteq \mathbb{R}^{n}\right\}<\infty.
 $$
\end{theorem}

 Furthermore,  the family of all bounded open sets in
the inequality (e) of Theorem \ref{th3} in some  situation can be
replaced by the family of all open balls.

\begin{theorem} \label{th4} Let  $\beta\in(0, n)$ and $\mu$ a nonnegative Radon measure on $\mathbb{R}^{1+n}_{+}.$ If
  $1<p<\min\{q, n/\beta\}$ or $1=p\leq q<\infty,$ then the following two conditions are equivalent:\\
(a)
$$
\|v(t^{2\alpha},x)\|_{L^{q}(\mathbb{R}_{+}^{1+n}, \mu)}\lesssim
\|v_{0}\|_{\dot{W}^{\beta,p}}, \forall v_{0}\in
\dot{W}^{\beta,p}(\mathbb{R}^{n}).
$$
 (b)
 $$
 \sup\limits_{x\in\mathbb{R}^{n}, r>0}
 \frac{\left(\mu\left(T(B(x,r))\right)\right)^{p/q}}{cap_{\dot{W}^{\beta,p}}\left(B(x,r)\right)}<\infty.
 $$
 But, this equivalence fails to hold  when $1<p=q< n/\beta.$
\end{theorem}
\begin{remark} We  plan to check in our future work that
whether or not the operator
$$
v_{0}\longrightarrow S_{t^{2\alpha}}v_{0}(x)
$$
 being compact from $\dot{W}^{\beta, p}(\mathbb{R}^{n})$  to
$L^{q}(\mathbb{R}^{1+n}_{+}, \mu)$ is equivalent to
 \begin{equation}\label{remeq1}
\lim\limits_{t\longrightarrow0}\frac{t^{p/q}}{c_{p}^{\beta}(\mu;
t)}=0,\ \hbox{ if}\ 1\leq p\leq n/\beta, p\leq q<\infty;
\end{equation}
 \begin{equation}\label{remeq2}
 \int_{0}^{\infty}\left(\frac{t^{p/q}}{c_{p}^{\beta}(\mu; t)}\right)^{q/(p-q)}\frac{d
 t}{t}<\infty\
\hbox{if}\ 0<q<p, 1< p\leq n/\beta;
\end{equation}
 or when $1<p<\min\{q, n/\beta\}$ or $1=p\leq q<\infty,$
 \begin{equation}\label{remeq3}
\lim\limits_{\delta\longrightarrow0}\sup\limits_{x\in\mathbb{R}^{n},
r\in(0,\delta)}
 \frac{\left(\mu\left(T(B(x,r))\right)\right)^{p/q}}{cap_{\dot{W}^{\beta,p}}\left(B(x,r)\right)}=0
\end{equation} and
 \begin{equation}\label{remeq4}
\lim\limits_{|x|\longrightarrow0}\sup\limits_{ r\in(0,1)}
 \frac{\left(\mu\left(T(B(x,r))\right)\right)^{p/q}}{cap_{\dot{W}^{\beta,p}}\left(B(x,r)\right)}=0.
 \end{equation}
 Similar results hold for the
embedding (\ref{equavi1}), see Maz'ya \cite[section 8.5, 8.6]{V.G.
Maz'ya 1} and \cite{V.G. Maz'ya 4} or Adams-Hedberg \cite[section
7.3]{D.R. Adams Hedberg}.
\end{remark}
\begin{remark} Since the inequality (b) in Theorem \ref{th4}
corresponds to the classical Carleson criterion for
$L^{p}(\mathbb{R}^{n})$ to be embedded in $L^{p}(\mathbb{R}^{1+n},
\mu)$ via Poisson's kernel (see for example Grafakos  \cite[ p. 539,
Theorem 7.37]{L. Grafakos}),
 we can refer to  the embeddings  in Theorems \ref{th2}-\ref{th4}  as the Carleson embeddings for the homogeneous  Sobolev spaces per the
 Cauchy problem for the $\alpha-$parabolic equation.
 \end{remark}

When $0<q<p=1$ we  obtain  necessary conditions for such embeddings.
\begin{theorem}\label{th5} Let  $\beta\in(0,n),$ $0<q<p=1$  and $\mu$ a nonnegative Radon  measure on $\mathbb{R}^{1+n}_{+}.$
 Then
 $(a)\Longrightarrow(b)\Longrightarrow(c)\Longrightarrow(d):$\\
 (a)
$$
\|v(t^{2\alpha},x)\|_{L^{q}(\mathbb{R}_{+}^{1+n}, \mu)}\lesssim
\|v_{0}\|_{\dot{W}^{\beta,1}}, \forall v_{0}\in
\dot{W}^{\beta,1}(\mathbb{R}^{n}).
$$
 (b)
 $$
 \|v(t^{2\alpha},x)\|_{L^{q,\infty}(\mathbb{R}^{1+n}_{+}, \mu)}\lesssim\|v_{0}\|_{\dot{W}^{\beta,1}(\mathbb{R}^{n})}, \
 \forall v_{0}\in \dot{W}^{\beta, 1}(\mathbb{R}^{n}).
 $$
 (c)
 $$
 \sup\left\{\frac{\left(\mu\left(T(O)\right)\right)^{1/q}}
 {cap_{\dot{W}^{\beta,1}}(O)}:\ \hbox{ open }\ O\subseteq \mathbb{R}^{n}\right\}<\infty.
 $$
(d)
$$
\|v(t^{2\alpha},x)\|_{L^{q,1}(\mathbb{R}^{1+n}_{+}, \mu)}\lesssim
\|v_{0}\|_{\dot{W}^{\beta,1}(\mathbb{R}^{n})}, \ \forall v_{0}\in
\dot{W}^{\beta, 1}(\mathbb{R}^{n}).
$$
 \end{theorem}

We can establish  the following  decay  of the  solutions  of
equation (\ref{a1}).
\begin{theorem} \label{th6}  If $v_{0}\in\dot{W}^{\beta,p}(\mathbb{R}^{n})$
 for $1\leq p<n/\beta$ and $\beta\in(0,n)$, then
$$
|v(t_{0}^{2\alpha},x_{0})|\lesssim
t_{0}^{p\beta-n}\|v_{0}\|_{\dot{W}^{\beta,p}(\mathbb{R}^{n})},\
\forall (t_{0},x_{0})\in \mathbb{R}^{1+n}_{+}.
$$
Equivalently
$$
|v(t_{0}, x_{0})|\lesssim
t_{0}^{\frac{p\beta-n}{2\alpha}}\|v_{0}\|_{\dot{W}^{\beta,p}(\mathbb{R}^{n})},\
\forall (t_{0},x_{0})\in \mathbb{R}^{1+n}_{+}.
$$

\end{theorem}

The special case $\alpha=p=1$ of Theorem \ref{th6} was  proved by
Xiao in \cite{J.Xiao 2}.

Working from $\mathbb{R}^{1+n}_{+}$ to $\mathbb{R}^{n},$ a trace
inequality can be derived from $\dot{W}^{\beta,p}(\mathbb{R}^{n}).$
\begin{theorem}\label{th7} Let $\beta\in (0,n),$  $1< p\leq q<\infty,$ $p<n/\beta$ and $\mu$ be a nonnegative Radon measure on $\mathbb{R}^{n}.$ Then
$$
\|f\|_{L^{q}(\mathbb{R}^{n},
\mu)}\lesssim\|f\|_{\dot{W}^{\beta,p}(\mathbb{R}^{n})},\ f\in
\dot{W}^{\beta,p}(\mathbb{R}^{n})
\Leftrightarrow\sup\limits_{\hbox{open}\ O\subseteq
\mathbb{R}^{n}}\frac{(\mu(O))^{p/q}}{cap_{\dot{W}^{\beta,p}(\mathbb{R}^{n})}(O)}<\infty.
$$
If $1=p\leq q<\infty,$ or $1<p<\min\{q,n/\beta\}$ then
$$
\|f\|_{L^{q}(\mathbb{R}^{n},
\mu)}\lesssim\|f\|_{\dot{W}^{\beta,1}(\mathbb{R}^{n})}, \ f\in
\dot{W}^{\beta,p}(\mathbb{R}^{n}) \Leftrightarrow\sup_{x\in
\mathbb{R}^{n},
r>0}\frac{(\mu(B(x,r)))^{p/q}}{cap_{\dot{W}^{\beta},p}(B(x,r))}
<\infty.
$$

\end{theorem}
Similarly to Theorem \ref{th1}, we obtain the following result which
covers Theorem \ref{th7}.
 \begin{theorem}\label{th8} Let $\beta\in (0,n),$  $1< p\leq q<\infty,$ $p<n/\beta$ and $\mu$ be a nonnegative  Radon measure on $\mathbb{R}^{n}$.
 Then the following statements are equivalent:\\
(a)
$$
\|f\|_{L^{q,p}(\mathbb{R}^{n}, \mu)}\lesssim\|f\|_{\dot{W}^{\beta,
p}}(\mathbb{R}^{n}), \ \ f\in \dot{W}^{\beta, p}(\mathbb{R}^{n}),
$$
 (b)
$$
\|f\|_{L^{q}(\mathbb{R}^{n},
\mu)}\lesssim\|f\|_{\dot{W}^{\beta,p}(\mathbb{R}^{n})}, \ \ f\in
\dot{W}^{\beta, p}(\mathbb{R}^{n}),
$$
 (c)
 $$
 \|f\|_{L^{q,\infty}(\mathbb{R}^{n}, \mu)}\lesssim\|f\|_{\dot{W}^{\beta, p}(\mathbb{R}^{n})}, \ \
f\in \dot{W}^{\beta, p}(\mathbb{R}^{n}),
$$
 (d)
 $$
 (\mu(O))^{p/q}\lesssim cap_{\dot{W}^{\beta, p}(\mathbb{R}^{n})}(O), \
 \hbox{open}\ O\subseteq \mathbb{R}^{n}.
 $$
 If $1=p\leq q<\infty,$  or $1<p<\min\{q,n/\beta\}$ then  all of them are equivalent to\\
 (e)
   $$
  \sup\limits_{r>0,x\in \mathbb{R}^{n}}\frac{(\mu(B(x,r)))^{p/q}}{cap_{\dot{W}^{\beta,p}(\mathbb{R}^{n})}(B(x,r))}<\infty.
 $$
If $0<q<p=1,$ then
$(b)\Longrightarrow(c)\Longrightarrow(d)\Longrightarrow(a).$
\end{theorem}
\begin{remark}
The equivalences $(b)\Longleftrightarrow(d)$ and
$(b)\Longleftrightarrow(e)$ in Theorem \ref{th8} can be verified
directly from Cascan-Ortega-Verbitsky \cite[Theorem 3.1 \&
3.2]{Cascan ortega verbitsky},
 Maz'ya \cite{V.G. Maz'ya 62}-\cite{V.G. Maz'ya 64} and \cite{V.G.
Maz'ya 1}, and Adams-Hedberg \cite[Theorem 7.2.2]{D.R. Adams
Hedberg}. Moreover, the  case $1=p\leq q<\infty$ of Theorem
 \ref{th8} was shown by Xiao in  \cite{J.Xiao 2}.
\end{remark}
Nishio-Yamada \cite{M. Nishio} gave a characterization for a
nonnegative Radon  measure $\mu$
 on $\mathbb{R}^{1+n}_{+}$
 to be  a Carleson type measure on $b_{\alpha}^{p}(\mathbb{R}^{1+n}_{+}),$
 which is
 called  (0,1)-type
Carleson measure and  means that
 $|\nabla_{(t,x)} u(t,x)|\in L^{p}(\mathbb{R}^{1+n}_{+},\mu),$  that is,
$$
\|\nabla_{(t,x)} u(t,x)\|_{L^{p}(\mathbb{R}^{1+n}_{+},
\mu)}\lesssim\|u(t,x)\|_{L^{p}(\mathbb{R}^{1+n}_{+})}, \ \forall
u\in b_{\alpha}^{p}(\mathbb{R}^{1+n}_{+}).
$$
We find a sufficient condition for a nonnegative Radon  measure
$\mu$
 on $\mathbb{R}^{1+n}_{+}$
 to be  a Carleson type measure on $b_{1/2}^{p}(\mathbb{R}^{1+n}_{+}).$

\begin{theorem}\label{th9} If  $\mu$ is a nonnegative Radon  measure  on $\mathbb{R}^{1+n}_{+}$ satisfying the property
$$
\sup\limits_{x\in \mathbb{R}^{n},
r>0}\frac{\left(\mu\left(T(B(x,r))\right)\right)^{p/q}}
 {cap_{\dot{W}^{1/2,p}}(B(x,r))}<\infty
 $$
for  $1\leq p<2n$ and $\frac{4pn+4p}{2n-p}\leq q<\infty,$ then $\mu$
is a (0,1)-type Carleson  measure  on
$b_{1/2}^{p_{1}}(\mathbb{R}_{+}^{1+n})$ for
$p_{1}=\frac{q(2n-p)}{2p(n+1)}-1.$
\end{theorem}

\begin{corollary} \label{co1} Let  $\beta\in(0,n),$ $1<p<\min\{q, n/\beta\}$ or $1=p\leq q<\infty,$ $\gamma\in (0,1],$ $\zeta>0$ and
$\zeta+n\gamma> n-p\beta.$ If
$d\mu_{\gamma,\zeta}(t,x)=t^{\zeta-1}|x|^{n(\beta-1)}dtdx,$ then
$$
\left(\mu_{\gamma,\zeta}(T(O))\right)^{\frac{n-p\beta}{\zeta+n\gamma}}\lesssim
cap_{\dot{W}^{\beta,p}}(O),\
 \hbox{open}\ O\subseteq\mathbb{R}^{n}.
  $$
Equivalently
 $$
 \|v(t^{2\alpha},x)\|_{L^{\frac{(\zeta+n\gamma)p}{n-p\beta}}(\mathbb{R}_{+}^{1+n}, \mu_{\gamma,\zeta})}
 \lesssim\|v_{0}\|_{\dot{W}^{\beta,p}},
\forall v_{0}\in \dot{W}^{\beta,p}.
$$
\end{corollary}
{\bf{Proof}.} This assertion follows from the case
$q=p(\zeta+n\gamma)/(n-p\beta)$ and $\mu=\mu_{\gamma,\zeta}$ of
Theorem \ref{th3}. $\Box$

The first inequality of Corollary \ref{co1} is the iso$-$capacitary
inequality (see Maz'ya \cite{V.G. Maz'ya 4} for more)
 and the second is its analytic form attached to the kernel $K_{t^{2\alpha}}^{\alpha}(x).$ Both of them   were firstly
 stated by Xiao in  \cite{J.Xiao 2} for $\alpha=p=1.$

\begin{corollary} \label{co2} Let $\alpha\in(0,1],$ $\beta\in(0, n),$ $1\leq p< n/\beta$ and $\gamma\in(-1, \infty).$ Then the following two
conditions hold:\\
(a)
$$
\left(\int_{\mathbb{R}^{1+n}_{+}}|v(t^{2\alpha},x)
|^{\frac{p(1+n+\gamma)}{n-p\beta}}t^{\gamma}dtdx\right)^{\frac{n-p\beta}{p(1+n+\gamma)}}
\lesssim\|v_{0}\|_{\dot{W}^{\beta,p}}, \forall v_{0}\in
\dot{W}^{\beta, p}(\mathbb{R}^{n}).
$$
(b)
$$
\sup\limits_{t>0}\left(\int_{\mathbb{R}^{n}}|v(t^{2\alpha},x)
|^{\frac{pn}{n-p\beta}}dx\right)^{\frac{n-p\beta}{pn}}\lesssim\|v_{0}\|_{\dot{W}^{\beta,p}},
 \forall v_{0}\in \dot{W}^{\beta, p}(\mathbb{R}^{n}).
 $$
\end{corollary}
{{\bf{Proof.}}} In Theorem \ref{th4} we take
$$
d\mu(t,x)=(1+\gamma)^{-1}t^{\gamma}dtdx, \
q=\frac{p(1+n+\gamma)}{n-p\beta},\ \gamma>-1,
$$
respective
$$
d\mu(t,x)=\delta_{t_{0}}(t)\otimes dx,\ q=\frac{pn}{n-p\beta},\
\gamma\longrightarrow-1,
$$
where $\delta_{t_{0}}(t)$ is the Dirac measure at $t_{0}>0,$ then an
application of the capacitary estimate of
 ball (see Maz'ya \cite[p. 356]{V.G. Maz'ya 1} for $p\in (1,n/\beta)$, Xiao \cite[p. 833]{J.Xiao 2} for $p=1$):
$$
cap_{\dot{W}^{\beta,p}}(B(x,r))\approx r^{n-p\beta}, x\in
\mathbb{R}^{n}, r>0,
$$
we can finish the proof. $\Box$

 According to Miao-Yuan-Zhang  \cite[ Proposition 2.1]{C. Miao}, the condition (a) of Corollary \ref{co2} amounts to that
  $\dot{W}^{\beta, p}(\mathbb{R}^{n})$  is embedded in the homogeneous Besov or Triebel-Lizorkin   space (see Triebel
   \cite{H. Triebel}  for more details about these spaces)
 $$
 \dot{B}^{-\frac{\gamma+1}{q}}_{q,q}(\mathbb{R}^{n})=\dot{F}^{-\frac{\gamma+1}{q}}_{q,q}(\mathbb{R}^{n}),\
  q=\frac{p(1+n+\gamma)}{n-p\beta}.
  $$
At the same time, the condition (b) of Corollary \ref{co2} can be
treated as extreme case of the condition (a) in Corollary \ref{co2}.

 The rest of this paper is organized as follows. In the next section,
 we  give six preliminary results: Lemma \ref{le1}$-$a strong-type inequality
 for the Hardy-littlewood maximal operator with respect to the
 variational capacity (whose new generalizations were made by Maz'ya \cite{Maz'ya 9} and  Costea-Maz'ya\cite{Costea Maz'ya}), Lemma \ref{le2}$-$an elementary Riesz integral  upper estimate of the kernel $K^{\alpha}_{t}(x)$, Lemma \ref{le3}$-$a lower
 estimate for the kernel $K^{\alpha}_{t}(x)$, Lemma \ref{le4}$-$four standard estimates involving capacity, measure and
 nontangential maximal functions,  Lemma \ref{le5}$-$an integral representation of fractional order  homogeneous Sobolev functions,  and
 Lemma \ref{le6}$-$a homogeneous version of the extension/restriction theorem.  In the third section, we prove our theorems and corollaries:
 Theorem \ref{th1}  is proved by using  Lemma \ref{le1}. Theorem \ref{th2} is showed from   applying  Lemmas \ref{le1}, \ref{le4}, \ref{le6}
   and the dyadically discrete
 forms of the left-hand integrals in $(a)-(b)$ of  Theorem \ref{th2}. Theorem \ref{th3} is demonstrated by using Lemmas \ref{le1} \& \ref{le4}.
   Theorem \ref{th4}
 is derived from  the equivalence established in Theorem \ref{th3}, Lemmas \ref{le2}-\ref{le5} and more delicate estimates for measures, functions and
 integrals under consideration.  Theorem \ref{th5} is verified  from  Lemmas \ref{le1} \& \ref{le4}.
  Theorem \ref{th6}  is obtained  by applying Theorems \ref{th3}
  \& \ref{th4}  and estimating the norm of Dirac measure on $\mathbb{R}^{1+n}_{+}.$  Theorem \ref{th9}  is established through comparing
   $1/2-$parabolic rectangles  in $\mathbb{R}^{1+n}_{+}$ with the tents of $n-$dimensional balls.

 \vspace{0.10in}
\section{{{Preliminary Lemmas}}}
\ \ This section contains six technical  results needed for proving
the main results of this paper. The first is the capacity
strong-type
 inequalities for $f\in \dot{W}^{\beta,p}(\mathbb{R}^{n})$ 
  and its Hardy-Littlewood maximal
  operator
$$
\mathcal
{M}f(x)=\sup\limits_{r>0}r^{-n}\int_{B(x,r)}|f(y)|\hbox{d}y, \ x\in
\mathbb{R}^{n}.
$$
\begin{lemma}\label{le1}
  The following three inequalities hold:\\
(a) If $\beta\in(0,n)$ and  $p\in[1, n/\beta],$  then, $ \forall
f\in \dot{W}^{\beta, p}(\mathbb{R}^{n}),$
$$
\int_{0}^{\infty}cap_{\dot{W}^{\beta,p}(\mathbb{R}^{n})}(\{x\in
\mathbb{R}^{n}: |f(x)|\geq \lambda\})d\lambda^{p} \lesssim
\|f\|_{\dot{W}^{\beta,p}(\mathbb{R}^{n})}^{p} ;
$$
 If $1\leq p<\infty,$ then, $ \forall
f\in \dot{W}^{1, p}(\mathbb{R}^{1+n}_{+}),$
$$
\int_{0}^{\infty}cap_{\dot{W}^{1,p}(\mathbb{R}^{1+n}_{+})}(\{(t,x)\in
\mathbb{R}^{1+n}_{+}: |f(t,x)|\geq\lambda\})d\lambda^{p} \lesssim
\|f\|_{\dot{W}^{1,p}(\mathbb{R}^{1+n}_{+})}^{p}.
$$
 (b) If $\beta\in(0,n)$ and  $p\in[1, n/\beta],$ then, $\forall
f\in \dot{W}^{\beta, p}(\mathbb{R}^{n}),$
$$
\int_{0}^{\infty}cap_{\dot{W}^{\beta,p}(\mathbb{R}^{n})}(\{x\in
\mathbb{R}^{n}: |\mathcal {M}f(x)|\geq \lambda\})d \lambda^{p}
\lesssim \|f\|_{\dot{W}^{\beta,p}(\mathbb{R}^{n})}^{p}.
$$
\end{lemma}
 {\bf{Proof.}} (a) {{\it{Part 1}}},  $f\in \dot{W}^{1, p}(\mathbb{R}^{1+n}_{n}):$ This assertion is due to
  Maz'ya \cite[Section 2.3.1]{V.G. Maz'ya 1} or his  another work \cite{V.G. Maz'ya 2}.
   {{\it{Part 2}}}, $f\in \dot{W}^{\beta, p}(\mathbb{R}^{n}):$
  {{{\em{Case 1}}},} $p\in (1, n/\beta):$ This case is due  to Maz'ya \cite[Proposition 4.1]{V.G. Maz'ya 3}
  or Maz'ya  \cite[p. 368 Theorem]{V.G. Maz'ya 1}.
 {{{{\em{Case 2}}},}} $p=1:$ This case is essentially prove by Wu \cite{Z. Wu} when $\beta\in(0,1)$ and
  Xiao \cite{J.Xiao 2} when $\beta \in(0, n)$.
   {{{{{\em{Case 3}}}},}} $p=n/\beta:$   It can be found in Maz'ya  \cite{Maz'ya} or Adams-Xiao \cite{D. R. Adams J. Xiao}.\\
 (b) If $f\in\dot{W}^{\beta, p}(\mathbb{R}^{n}):$ We divide the proof into three cases.\\
 {{{{\em{Case 1}},}}} $p=1:$ It is due to  Xiao \cite{J.Xiao 2}.
 {{{{\em{Case 2}}}},} $p=n/\beta:$  This is  proved by Adams-Xiao \cite{D. R. Adams J. Xiao}.
{{{{{\em{Case 3}}}}},} $p\in(1, n/\beta):$ It follows from Maz'ya
\cite[p. 347, Theorem 2]{V.G. Maz'ya 1} or his earlier work
  \cite{V.G. Maz'ya Havin}  that for $1<p<n/\beta,$ $f\in  \dot{W}^{\beta, p}(\mathbb{R}^{n})$ if and only if
$$
f=(-\triangle)^{-\beta/2}g=I_{\beta}\ast g(x)\ \hbox{ and}\ \
 \|f\|_{\dot{W}^{\beta,p}}=\|g\|_{L^{p}},
 $$
 for some $g\in L^{p}(\mathbb{R}^{n})$, where
 $$
 I_{\beta}\ast g(x)=\frac{1}{\gamma_{\beta}}\int_{\mathbb{R}^{n}}\frac{g(y)}{|x-y|^{n-\beta}}dy
 $$
  with  $\gamma_{\beta}=\pi^{n/2}2^{\beta}\Gamma(\beta/2)/\Gamma(\frac{n-\beta}{2}).$ Then for fixed
  $f\in \dot{W}^{\beta, p}(\mathbb{R}^{n}),$ and  $g\in L^{p}(\mathbb{R}^{n})$ with $f(x)=I_{\beta}\ast g(x),$
  according to  R. Johnson \cite[p. 33, Proof of Theorem 1.9]{R. Johnson}, we have
 $$
 \mathcal{M}(I_{\beta}\ast g)\leq I_{\beta}\ast(\mathcal{M}{g})
 $$
and
 $$ M_{\lambda}(\mathcal {M}f(x))\subseteq
M_{\lambda}(I_{\beta}\ast(\mathcal{M}( g)).
$$
 It follows from Maximal Theorem  Stein \cite[p. 13, Theorem 1]{E.M. Stein} that
 $$
 \mathcal{M}( g)\in L^{p}(\mathbb{R}^{n})\ \hbox{ and}\
\|\mathcal{M}(g)\|_{p}\lesssim\|g\|_{p}.
$$
  Thus (a) implies (b). $\Box$

\begin{lemma}\label{le2}
If  $\alpha\in(0,1],$ $\beta\in (0,n)$ and
$(t,x)\in\mathbb{R}_{+}^{1+n},$ then
$$
\int_{\mathbb{R}^{n}}K_{t^{2\alpha}}^{\alpha}(y)|y-x|^{\beta-n}dy
\lesssim (t^{2}+|x|^{2})^{\frac{\beta-n}{2}}.
$$
\end{lemma}
 {\bf{Proof.}} By  Miao-Yuan-Zhang \cite{C. Miao},  Nishio-Shimomura-Suzuki \cite{M. Nishio K. Shimomura N. Suzuki}
  or Nishio-Yamada \cite{M. Nishio}, we have the following point-wise estimate
   \begin{equation}
|K_{t}^{\alpha}(x)|\leq C\frac{t}{(t^{1/2\alpha}+|x|)^{n+2\alpha}},\
\forall (t,x)\in \mathbb{R}^{1+n}_{+}.\label{a3}
 \end{equation}
 So, it  suffices to verify
$$
J(t,x):=\int_{\mathbb{R}^{n}}t^{2\alpha}(t+|y|)^{-n-2\alpha}|y-x|^{\beta-n}dy\lesssim
(t^{2}+|x|^{2})^{\frac{\beta-n}{2}}.
$$
 Changing variables: $x\longrightarrow tx,$ $y\longrightarrow ty,$ we see the previous estimate is equivalent to the following one:
$$
J(1,x)\lesssim (1+|x|^{2})^{\frac{\beta-n}{2}}.
$$
Since $J(1,0)\lesssim 1$ we  may assume that $|x|>0.$ Then
  \begin{equation}
 J(1,x)\nonumber\lesssim
 \left(\int_{B(x, |x|/2)}+\int_{\mathbb{R}^{n}\backslash B(x, |x|/2)}\right)
 \frac{1}{(1+|y|)^{n+2\alpha}|y-x|^{n-\beta}}dy=I_{1}(x)+I_{2}(x).\nonumber
  \end{equation}
 Since  $|x-y|\leq|x|/2$ implies that $|y|\approx |x|,$ we have
 \begin{eqnarray*}
 I_{1}(x)\nonumber
 &=&\int_{B(x,|x|/2)}\frac{1}{(1+|y|)^{n+2\alpha}|y-x|^{n-\beta}}dy\\ \nonumber
&\lesssim&(1+|x|)^{-n-2\alpha}\int_{B(x,|x|/2)}\frac{1}{|y-x|^{n-\beta}}dy\\\nonumber
&\lesssim&(1+|x|)^{-n-2\alpha}\int_{0}^{|x|/2}s^{\beta-1}ds\\\nonumber
&\lesssim& |x|^{\beta}(1+|x|)^{-n-2\alpha}\\\nonumber
 &\lesssim&(1+|x|)^{\beta-n},\nonumber
  \end{eqnarray*}
with the last inequality using the fact $1\leq (1+|x|)^{2\alpha}.$
If $|x-y|> |x|/2,$ then
 \begin{eqnarray*} I_{2}(x)\nonumber
&=&\int_{\mathbb{R}^{n}\backslash B(x,
|x|/2)}\frac{1}{(1+|y|)^{n+2\alpha}|x-y|^{n-\beta}}dy\\\nonumber
&\lesssim&|x|^{\beta-n}\int_{\mathbb{R}^{n}\backslash
B(x,|x|/2)}\frac{1}{(1+|y|)^{n+2\alpha}}dy\\\nonumber
&\lesssim&|x|^{\beta-n},\nonumber
\end{eqnarray*}
with the last inequalities using the fact
$\frac{1}{(1+|y|)^{n+2\alpha}}\in L^{1}(\mathbb{R}^{n}).$ Since
$|x-y|>|x|/2$ implies
 $|y|<3|x-y|,$
 $$
 I_{2}(x)\lesssim\int_{\mathbb{R}^{n}\backslash B(x, |x|/2)}\frac{1}{(1+|y|)^{n+2\alpha}|y|^{n-\beta}}dy\lesssim 1.
 $$
 Thus $I_{2}\lesssim(1+|x|)^{\beta-n}$ and $J(1,x)\lesssim(1+|x|^{2})^{\frac{\beta-n}{2}}.$ $\Box$

\begin{lemma}\label{le3}\cite{M. Nishio} For $\alpha\in(0,1],$
there are positive constants $\sigma$ and $C$ such that
$$
\inf\{|K_{t}^{\alpha}(x)|: |x|\leq \sigma
t^{\frac{1}{2\alpha}}\}\geq Ct^{-\frac{n}{2\alpha}},
$$
where $\sigma$ and $C$ depend only on $n, \alpha.$
 \end{lemma}

 \begin{lemma} \label{le4} Let $\alpha\in(0,1]$ and  $ \beta\in(0,n).$ Given $f\in \dot{W}^{\beta, p}(\mathbb{R}^{n}),$ $\lambda>0,$ and a
  nonnegative measure $\mu$ on $\mathbb{R}_{+}^{1+n},$ let
 $$
 E_{\lambda}^{\alpha,\beta}(f)=\{(t,x)\in \mathbb{R}_{+}^{1+n}: |S_{\alpha}(t^{2\alpha})f(x)|>\lambda\}
 $$
and
$$
O_{\lambda}^{\alpha,\beta}(f)=\{y\in \mathbb{R}^{n}:
\sup\limits_{|y-x|<t}|S_{\alpha}(t^{2\alpha})f(x)|>\lambda\}.
$$
Then the following four statements are true:\\
(a) For any natural number $k,$
$$
\mu\left(E_{\lambda}^{\alpha,\beta}(f)\cap T(B(0,k))\right)\leq
\mu\left(T(O_{\lambda}^{\alpha,\beta}(f)\cap B(0,k))\right).
$$
(b) For any natural number $k,$
$$
cap_{\dot{W}^{\beta,p}(\mathbb{R}^{n})}\left(O_{\lambda}^{\alpha,\beta}(f)\cap
B(0,k)\right)\geq c_{p}^{\beta} \left(\mu;
\mu\left(T(O_{\lambda}^{\alpha,\beta}(f)\cap B(0,k))\right)\right).
$$
(c) There exists a dimensional constant $\theta_{1}>0$ such that
$$
\sup\limits_{|y-x|<t}|S_{\alpha}(t^{2\alpha})f(y)|\leq
\theta_{1}\mathcal {M}f(x),  x\in \mathbb{R}^{n}.
$$
(d) There exists a dimensional constant $\theta_{2}>0$ such that
$$
(t,x)\in T(O)\Longrightarrow S_{\alpha}(t^{2\alpha})|f|(x)\geq
\theta_{2},
$$
where  $O$ is a bounded open set contained in $\hbox{Int}(\{x\in
\mathbb{R}^{n}: f(x)\geq 1 \}).$
  \end{lemma}
 {\bf{Proof.}} (a) Since  $\sup\limits_{|y-x|<t}|S_{\alpha}(t^{2\alpha})f(x)|$ is   lower semicontinuous on $\mathbb{R}^{n},$
$O_{\lambda}^{\alpha,\beta}(f)$ is an open subset of
$\mathbb{R}^{n}.$ By the definition of
$E_{\lambda}^{\alpha,\beta}(f)$ and $O_{\lambda}^{\alpha,\beta}(f),$
we have
$$
E_{\lambda}^{\alpha,\beta}(f)\subseteq
T(O_{\lambda}^{\alpha,\beta}(f))\ \hbox{and} \
\mu(E_{\lambda}^{\alpha,\beta}(f))\leq
T(\mu(O_{\lambda}^{\alpha,\beta}(f))).
$$
 Then
$$
\mu\left(E_{\lambda}^{\alpha,\beta}(f)\cap T(B(0,k))\right)\leq
\mu(T(O_{\lambda}^{\alpha,\beta}(f) \cap
T(B(0,k))))=\mu\left(T(O_{\lambda}^{\alpha,\beta}(f)\cap
B(0,k))\right).
$$
(b) It follows from the definition of $c_{p}^{\beta}(\mu; t).$ \\
(c) By (\ref{a3}), we have
$$
|S_{\alpha}(t^{2\alpha})f(x)|=|K_{t^{2\alpha}}^{\alpha}(x)\ast
f(x)|\leq \int_{\mathbb{R}^{n}}\frac{Ct^{2\alpha}}
{(t+|x-y|)^{n+2\alpha}}|f(y)|\hbox{d}y:=H_{t}(x)\ast|f(x)|.
$$
Thus
$$
\sup\limits_{|y-x|<t}|S_{\alpha}(t^{2\alpha})f(y)|\leq
\sup\limits_{|y-x|<t}H_{t}(y)\ast|f(y)|\leq
\theta_{1}\mathcal{M}f(x).
$$
 The last inequality follows from Stein \cite[p. 57, Proposition]{E.M. Stein}.\\
 (d) For any $(t,x)\in T(O),$ we have
$$
 B(x,t)\subseteq O\subseteq \hbox{Int}(\{x:f(x)\geq1\}).
 $$
It follows from  Lemma \ref{le3}  that there exist $\sigma$ and $C$
which are only depending on $n$ and $\alpha$ such that
$$
\inf\{K_{t}^{\alpha}(x): |x|\leq \sigma t^{\frac{1}{2\alpha}}\}\geq
Ct^{-\frac{n}{2\alpha}}.
$$
 Then
\begin{eqnarray*} S_{\alpha}(t^{2\alpha})|f|(x)\nonumber&=&\int_{\mathbb{R}^{n}}K_{t^{2\alpha}
 }^{\alpha}(x-y)|f|(y)\hbox{d}y\\\nonumber
&\geq& Ct^{-n}\int_{B(x,\sigma t)\cap \hbox{Int}(\{x:
 f(x)\geq
 1\})}|f|(y)\hbox{d}y.\nonumber
\end{eqnarray*}
If $\sigma>1,$ then
$$
B(x, \sigma t)\cap \hbox{Int}(\{x:
 f(x)\geq
 1\})\supseteq B(x, t)\cap \hbox{Int}(\{x:
 f(x)\geq
 1\})=B(x,t);
 $$
  if $\sigma \leq 1$ then
   $$
   B(x, \sigma t)\cap \hbox{Int}(\{x: f(x)\geq 1\})=B(x,\sigma t).
 $$
 Thus $S_{\alpha}(t^{2\alpha})|f|(x)\geq \theta_{2}$ for some  dimensional constant $\theta_{2}>0.$  $\Box$

Using
$f(x)=(-\triangle)^{-\beta/2}\left((-\triangle)^{\beta/2}f(x)\right)$
and the definition of Riesz potentials
 We can easily derive an integral representation of homogeneous Sobolev functions.
 \begin{lemma}\label{le5}   \cite{D.R. Adams Hedberg}   Let $p\in (1, n/\beta),$  $0<\beta<n$ and $f\in \dot{W}^{\beta,p}(\mathbb{R}^{n}).$
 Then
$$
f(x)=\frac{1}{\gamma_{\beta}}\int_{\mathbb{R}^{n}}\frac{(-\triangle)^{\beta/2}f(y)}{|y-x|^{n-\beta}}\hbox{d}y,
$$
where
$\gamma_{\beta}=\pi^{n/2}2^{\beta}\Gamma(\beta/2)/\Gamma(\frac{n-\beta}{2}).$
\end{lemma}

 The following result is a special case of Adamas \cite[Theorem 5.2]{D.R. Adams 4} or Adams-Xiao \cite[Theorem A]{D. R. Adams J. Xiao}.

 \begin{lemma}  \label{le6} Let $\beta\in (0,n).$ Then there are a linear extension operator
  $$
  \mathcal{E}: \dot{W}^{\beta, n/\beta}(\mathbb{R}^{n})\longrightarrow  \dot{\mathcal{L}}_{2\beta}^{n/\beta}(\mathbb{R}^{2n})
  $$
 and a linear restriction operator
 $$
 \mathcal{R}: \dot{\mathcal{L}}_{2\beta}^{n/\beta}(\mathbb{R}^{2n}) \longrightarrow \dot{W}^{\beta,n/\beta}(\mathbb{R}^{n})
 $$
  such that $\mathcal{R}\mathcal{E}$ is the identity, and \\
 (a)
  $$
  \|\mathcal{E}f\|_{\dot{\mathcal{L}}_{2\beta}^{n/\beta}(\mathbb{R}^{2n})}\lesssim \|f\|_{\dot{W}^{\beta, n/\beta}(\mathbb{R}^{n})},
  \ \forall f\in  \dot{W}^{\beta, n/\beta}(\mathbb{R}^{n});
 $$
(b)
$$
\|\mathcal{R}g\|_{\dot{W}^{\beta,
n/\beta}(\mathbb{R}^{n})}\lesssim\|g\|_{\dot{\mathcal{L}}_{2\beta}^{n/\beta}(\mathbb{R}^{2n})},\
\ \forall g\in
\dot{\mathcal{L}}_{2\beta}^{n/\beta}(\mathbb{R}^{2n}).
$$
 \end{lemma}

\vspace{0.1in}
\section{{{Proofs of Main Results}}}

\subsection{Proof of Theorem \ref{th1}} Assume that
$1\leq p\leq q<\infty.$ In what follows, for $\lambda>0$ and $ u\in
\dot{W}^{1,p}(\mathbb{R}^{1+n}_{+})\cap
b^{p}_{\alpha}(\mathbb{R}^{1+n}_{+}), $
 let
$$
M_{\lambda}(u)=\{(t,x)\in \mathbb{R}^{1+n}_{+}: |u(t,x)|\geq\lambda
\}.
$$
$(a)\Longrightarrow(b)\Longrightarrow(c).$ Since
$0<\lambda_{1}<\lambda_{2}$ implies $\mu(M_{\lambda_{2}}(u))\leq
\mu(M_{\lambda_{1}}(u)),$ we can conclude
$$
q\mu(M_{\lambda}(u))\lambda^{q-1}\leq\frac{d}{d\lambda}\left(\int_{0}^{\lambda}
(\mu(M_{s}(u)))^{p/q}ds^{p}\right)^{q/p}.
$$
 This implies
$$
\left(s^{q}\mu(M_{s}(u))\right)^{p/q}\leq\left(q\int_{0}^{\infty}
\mu(M_{\lambda}(u))\lambda^{q-1}d\lambda\right)^{p/q}\leq\int_{0}^{\infty}(\mu(M_{\lambda}(u)))^{p/q}d\lambda^{p},
\ \ s>0,
$$
and obtains  the desired implications.\\ $(c)\Longrightarrow(d).$
Let (c) be true. For an given  open  set $O\subseteq
\mathbb{R}^{1+n}_{+},$ and any function $u\in
\dot{W}^{1,p}(\mathbb{R}^{1+n}_{+})\cap
b^{p}_{\alpha}(\mathbb{R}^{1+n}_{+})$ with
 $$
 O\subseteq\hbox{Int}(\{(t,x)\in \mathbb{R}^{1+n}_{+}: u(t,x)\geq
1\}),
$$
we have $ \mu(O)\leq\mu(M_{1}(u))\lesssim\|u\|_{\dot{W}^{1,p}}^{q}.
$
 This derives  (d).\\
$(d)\Longrightarrow(a).$  If (d) is true, then for  $u\in
\dot{W}^{1,p}(\mathbb{R}^{1+n}_{+}),$ $k\in\mathbb{N}$ and
$B(0,k)\subseteq\mathbb{R}^{n},$  Lemma \ref{le1} (a) implies
\begin{eqnarray*}
&&\int\limits_{0}^{\infty}(\mu(M_{\lambda}(u)\cap ((0,k)\times
B(0,k))))^{p/q}d\lambda^{p}\\
 &\lesssim&\int\limits_{0}^{\infty}cap_{\dot{W}^{1,p}(\mathbb{R}^{1+n}_{+})}(M_{\lambda}(u)\cap
((0,k)\times B(0,k))))d\lambda^{p}\\
&\lesssim&\int\limits_{0}^{\infty}cap_{\dot{W}^{1,p}(\mathbb{R}^{1+n}_{+})}(M_{\lambda}(u))d\lambda^{p}\lesssim
\|u\|^{p}_{\dot{W}^{1,p}(\mathbb{R}^{1+n}_{+})}.
\end{eqnarray*}
 Letting $k\longrightarrow\infty$ we see (a) hold. When $0<q<p=1,$ the implications are obviously. $\Box$

 \vspace{0.1in}
\subsection{Proof of Theorem \ref{th2}}
 Let $0<q<p.$ Then we finish the proof in two steps.

{{\it{Part 3.2.1}}}: $(b)\Longrightarrow(a).$
 If
 $$
 I_{p,q}(\mu)=\int\limits_{0}^{\infty}\left(\frac{t^{p/q}}
 {c_{p}^{\beta}(\mu; t)}\right)^{\frac{q}{p-q}}\frac{d t}{t}<\infty,
 $$
 then for each $v_{0}\in\dot{W}^{\beta,p}(\mathbb{R}^{n}),$  each  $j=0,\pm1,\pm2,\cdots $ and each natural number $k,$
  Lemma \ref{le4} (c) implies
 $$
 cap_{\dot{W}^{\beta,p}(\mathbb{R}^{n})}(O_{2^{j}}(v_{0})\cap B(0,k))\leq
 cap_{\dot{W}^{\beta,p}(\mathbb{R}^{n})}(\{x\in \mathbb{R}^{n}: \theta_{1}\mathcal{M}v_{0}(x)>2^{j}\}\cap B(0,k)).
 $$
 Let $\mu_{j,k}(v_{0})=\mu(T(O_{2^{j}}(v_{0})\cap B(0,k)),$ and
 $$
 S_{p,q,k}(\mu; v_{0})=\sum_{j=-\infty}^{\infty}\frac{(\mu_{j,k}(v_{0})-\mu_{j+1,k}(v_{0}))^{\frac{p}{p-q}}}
 {\left(cap_{\dot{W}^{\beta,p}(\mathbb{R}^{n})}(O_{2^{j}}(v_{0})\cap B(0,k))\right)^{\frac{q}{p-q}}}.
 $$
 Lemma \ref{le4} (b) implies that
  \begin{eqnarray*}
  (S_{p,q,k}(\mu; v_{0}))^{\frac{p-q}{p}}\nonumber&=&
\left(\sum_{j=-\infty}^{\infty}\frac{(\mu_{j,k}(v_{0})-\mu_{j+1,k}(v_{0}))^{\frac{p}{p-q}}}
 {\left(cap_{\dot{W}^{\beta,p}(\mathbb{R}^{n})}(O_{2^{j}}(v_{0})\cap
 B(0,k))\right)^{\frac{q}{p-q}}}\right)^{\frac{p-q}{p}}\\\nonumber
&\lesssim&
\left(\sum_{j=-\infty}^{\infty}\frac{(\mu_{j,k}(v_{0})-\mu_{j+1,k}(v_{0}))^{\frac{p}{p-q}}}
 {\left(c_{p}^{\beta}(\mu; \mu_{j,k}(v_{0}))\right)^{\frac{q}{p-q}}}\right)^{\frac{p-q}{p}}\\\nonumber
&\lesssim&
\left(\sum_{j=-\infty}^{\infty}\frac{((\mu_{j,k}(v_{0}))^{\frac{p}{p-q}}-(\mu_{j+1,k}(v_{0}))^{\frac{p}{p-q}}}
 {\left(c_{p}^{\beta}(\mu; \mu_{j,k}(v_{0}))\right)^{\frac{q}{p-q}}}\right)^{\frac{p-q}{p}}\\\nonumber
&\lesssim& \left(\int\limits_{0}^{\infty}\frac{ds^{\frac{p}{p-q}}}
 {\left(c_{p}^{\beta}(\mu; s)\right)^{\frac{q}{p-q}}}\right)^{\frac{p-q}{p}}\\\nonumber
&\approx& (I_{p,q}(\mu))^{\frac{p-q}{p}}.
  \end{eqnarray*}
On the other hand,  using H$\ddot{o}$lder's inequality and Lemmas
\ref{le1} (b) and \ref{le4} (b)--(c), we have
\begin{eqnarray*}
&&\int_{T(B(0,k))}|v(t^{2\alpha}, x)|^{q}\hbox{d}\mu(t,x)\\
 &=&\int_{0}^{\infty}\mu\left(E_{\lambda}^{\alpha,\beta}(v_{0})\cap
T(B(0,k))\right)\hbox{d}\lambda^{q}\nonumber\\
 &\lesssim&\sum\limits_{j=-\infty}^{\infty}(\mu_{j,k}(v_{0})-\mu_{j+1,k}(v_{0}))2^{jq}\\\nonumber
&\lesssim&(S_{p,q,k}(\mu; v_{0}))^{\frac{p-q}{p}}
\left(\sum\limits_{j=-\infty}^{\infty}2^{jp}cap_{\dot{W}^{\beta,p}
(\mathbb{R}^{n})}(O_{2^{j}}(v_{0})\cap B(0,k))\right)^{q/p}\\
\nonumber
 & \lesssim&(S_{p,q,k}(\mu;
v_{0}))^{\frac{p-q}{p}}\!\left(\!\sum\limits_{j=-\infty}^{\infty}\!2^{jp}cap_{\dot{W}^{\beta,p}(\mathbb{R}^{n})}(\{x\in
\mathbb{R}^{n}: \theta_{1}\mathcal{M}v_{0}(x)>2^{j}\}\cap
B(0,k))\!\!\right)^{\frac{q}{p}}\\
 \nonumber
  &\lesssim&(S_{p,q,k}(\mu;
v_{0}))^{\frac{p-q}{p}}\left(\int\limits_{0}^{\infty}cap_{\dot{W}^{\beta,p}(\mathbb{R}^{n})}(\{x\in
\mathbb{R}^{n}:
\theta_{1}\mathcal{M}v_{0}(x)>\lambda\})\hbox{d}\lambda^{p}\right)^{q/p}\\\nonumber
&\lesssim& (S_{p,q,k}(\mu;
v_{0}))^{\frac{p-q}{p}}\|v_{0}\|^{q}_{\dot{W}^{\beta,
p}(\mathbb{R}^{n})}.
 \end{eqnarray*}
 Hence
$$
\left(\int_{T(B(0, k))}|v(t^{2\alpha},
x)|^{q}\hbox{d}\mu(t,x)\right)^{1/q} \lesssim
(I_{p,q}(\mu))^{\frac{p-q}{pq}}\|v_{0}\|_{\dot{W}^{\beta,
p}(\mathbb{R}^{n})}.
$$
 Letting $k\longrightarrow \infty$ in the left side of the above  estimate, we have
$$
\left(\int_{\mathbb{R}^{1+n}_{+}}|v(t^{2\alpha},
x)|^{q}\hbox{d}\mu(t,x)\right)^{1/q} \lesssim
(I_{p,q}(\mu))^{\frac{p-q}{pq}}\|v_{0}\|_{\dot{W}^{\beta,
p}(\mathbb{R}^{n})}.
$$

 {{\it{Part 3.2.2}}}: $(a)\Longrightarrow(b).$

If (a) is true, then
 $$
 J_{p, q}(\mu)=\sup\limits_{v_{0}\in \dot{W}^{\beta, p}(\mathbb{R}^{n}),\|v_{0}\|_{\dot{W}^{\beta,
p}(\mathbb{R}^{n})}>0}
 \frac{\left(\int_{\mathbb{R}_{+}^{1+n}}|v(t^{2\alpha},x)|^{q}\hbox{d}\mu(t,x)\right)^{1/q}}
 {\|v_{0}\|_{\dot{W}^{\beta,
p}(\mathbb{R}^{n})}}<\infty.
$$
 Thus for each $v_{0}\in \dot{W}^{\beta, p}(\mathbb{R}^{n})$, with $\|v_{0}\|_{\dot{W}^{\beta, p}(\mathbb{R}^{n})}>0,$
 we have
 $$
 \left(\int_{\mathbb{R}_{+}^{1+n}}|v(t^{2\alpha},x)|^{q}\hbox{d}\mu(t,x)\right)^{1/q}
 \leq J_{p,q}(\mu)\|v_{0}\|_{\dot{W}^{\beta,
p}(\mathbb{R}^{n})}.
$$
Since $\mu(E_{\lambda}^{\alpha,\beta}(v_{0}))$ is nonincreasing in
$\lambda,$  we have
 \begin{equation}
 \sup\limits_{\lambda>0}\lambda\left(\mu(E_{\lambda}^{\alpha,\beta}(v_{0}))^{1/q}\right)
 \lesssim J_{p,q}(\mu)\|v_{0}\|_{\dot{W}^{\beta,
p}(\mathbb{R}^{n})}.\label{a4}
 \end{equation}
 For   fixed positive $v_{0}\in \dot{W}^{\beta,p}(\mathbb{R}^{n}),$ and a bounded  open set
  $O\subseteq\hbox{Int}(\{x\in \mathbb{R}^{n}: v_{0}(x)\geq1\}),$ then (\ref{a4}) and Lemma \ref{le4} (d) imply that
 $$
 \mu(T(O))\leq \mu(E_{\frac{\theta_{2}}{2}}^{\alpha,\beta}(v_{0}))\lesssim (J_{p,q}(\mu))^{q}\|v_{0}\|^{q}_{\dot{W}^{\beta,p}(\mathbb{R}^{n})}.
 $$
 This along with the definition of $cap_{\dot{W}^{\beta,p}(\mathbb{R}^{n})}(\cdot)$ give
\begin{equation}
 \mu(T(O)) \lesssim
 (J_{p,q}(\mu))^{q}\left(cap_{\dot{W}^{\beta,p}(\mathbb{R}^{n})}(O)\right)^{q/p}.\label{a5}
 \end{equation}
 It follows from (\ref{a5}) and the definition of $c_{p}^{\beta}(\mu; t)$ that for $0<t<\infty,$  $c_{p}^{\beta}(\mu; t)>0.$
 The definition of $c_{p}^{\alpha}(\mu; t)$ implies that for every  integer $j$ there exists a bounded open set $O_{j}\subseteq
 \mathbb{R}^{n}$ such that
 $$
 cap_{\dot{W}^{\beta,p}(\mathbb{R}^{n})}(O_{j})\leq 2 \ c_{p}^{\beta}(\mu; 2^{j})\ \hbox{and}\ \mu(T(O_{j}))>2^{j}.
 $$
We divide the following proof into two cases.\\
{{\em{Case 1}},} $p\in (1, n/\beta):$\\
 It follows from Maz'ya \cite{V.G. Maz'ya 1} that
  $$
  cap_{\dot{W}^{\beta,p}(\mathbb{R}^{n})}(S)\approx\inf\left\{\|g\|^{p}:
 g\in L^{p}(\mathbb{R}^{n}), g\geq 0, S\subseteq
 \hbox{Int}(\{x\in\mathbb{R}^{n}: I_{\beta}\ast g(x)\geq 1\})\right\}.
 $$
By this equivalent  definition  we can find $g_{j}(x)\in L^{
p}(\mathbb{R}^{n})$ such that
 $$
 g_{j}\geq0, I_{\beta}\ast g_{j}(x)\geq1,\forall x\in O_{j} \ \hbox{and}\
  \|g_{j}\|_{L^{p}}^{p}\leq2\ cap_{\dot{W}^{\beta,p}(\mathbb{R}^{n})}(O_{j})\leq 4\ c_{p}^{\beta}(\mu; 2^{j}).
  $$
Given integers $i, k$ with $i<k,$ define
 $$
 g_{i,k}=\sup_{i\leq j\leq
k}\left(\frac{2^{j}}{c_{p}^{\beta}(\mu;2^{j})}\right)^{\frac{1}{p-q}}g_{j}.
$$
Since $L^{p}(\mathbb{R}^{n})$ is a lattice, we can conclude that
$g_{i,k}\in L^{p}(\mathbb{R}^{n})$ and
$$
\|g_{i,k}\|_{L^{p}}^{p}\leq
\sum\limits_{j=i}^{k}\left(\frac{2^{j}}{c_{p}^{\beta}(\mu;2^{j})}\right)^{\frac{p}{p-q}}
\|g_{j}\|_{L^{p}}^{p}\lesssim
\sum\limits_{j=i}^{k}\left(\frac{2^{j}}{c_{p}^{\beta}(\mu;2^{j})}\right)^{\frac{p}{p-q}}c_{p}^{\beta}(\mu;
2^{j}).
$$ Note that for $i\leq j\leq k,$
$$
x\in O_{j}\Longrightarrow I_{\beta}\ast
g_{i,k}(x)\geq\left(\frac{2^{j}}{c_{p}^{\beta}(\mu;2^{j})}\right)^{\frac{1}{p-q}}.
$$
It follows from Lemma \ref{le4} (d) that there exists  a dimensional
constant $\theta_{2}$ such that
$$
(t,x)\in T(O_{j})\Longrightarrow
S_{\alpha}(t^{2\alpha})|I_{\beta}\ast g_{i,k}(x)|(x)
\geq\left(\frac{2^{j}}{c_{p}^{\beta}(\mu;2^{j})}\right)^{\frac{1}{p-q}}\theta_{2}.
$$
This gives
$$
2^{j}<\mu\left(T(O_{j})\right)\leq\mu\left(E^{\alpha,\beta}_{\left(\frac{2^{j}}{c_{p}^{\beta}
(\mu;2^{j})}\right)^{\frac{1}{p-q}}\left(\frac{\theta_{2}}{2}\right)}(I_{\beta}\ast
g_{i,k}(x))\right).
$$
 Thus
  \begin{eqnarray*}
\left(J_{p,q}(\mu)\|g_{i,k}\|_{L^{p}}\right)^{q}\nonumber &\gtrsim&
\int_{\mathbb{R}_{+}^{1+n}}|S_{\alpha}(t^{2\alpha})(I_{\beta}\ast
g_{i,k}(x))|^{q}\hbox{d}\mu(t,x)\\\nonumber &\approx&
\int_{0}^{\infty}\left(\inf\{\lambda:
\mu\left(E_{\lambda}^{\alpha,\beta}(I_{\beta}\ast
g_{i,k}(x))\right)\leq s\}\right)^{q}\hbox{d}s\\\nonumber &\gtrsim&
\sum\limits_{j=i}^{k}\left(\inf\{\lambda:
\mu\left(E_{\lambda}^{\alpha,\beta}(I_{\beta}\ast
g_{i,k}(x))\right)\leq 2^{j}\}\right)^{q}2^{j}\\\nonumber
&\gtrsim&\sum\limits_{j=i}^{k}
\left(\frac{2^{j}}{c_{p}^{\beta}(\mu;2^{j})}\right)^{\frac{q}{p-q}}2^{j}\\\nonumber
&\gtrsim&\left(\frac{\sum_{j=i}^{k}\left(\frac{2^{j}}{c_{p}^{\beta}(\mu;2^{j})}\right)^{\frac{q}{p-q}}2^{j}}
{\left(\sum_{j=i}^{k}\left(\frac{2^{j}}{c_{p}^{\beta}(\mu;2^{j})}\right)^{\frac{q}{p-q}}c_{p}^{\beta}(\mu;
2^{j})\right)^{\frac{q}{p}}}\right)\|g_{i,k}\|_{L^{p}}^{q}\\\nonumber
&\approx&\left(\sum\limits_{j=i}^{k}\frac{2^{\frac{jp}{p-q}}}{\left(c_{p}^{\beta}(\mu;
2^{j})\right)^{\frac{q}{p-q}}}\right)^{\frac{p-q}{p}}\|g_{i,k}\|_{L^{p}}^{q}.\nonumber
 \end{eqnarray*}
 This tells us
 $$
 \sum\limits_{j=i}^{k}\frac{2^{\frac{jp}{p-q}}}{\left(c_{p}^{\beta}(\mu;
2^{j})\right)^{\frac{q}{p-q}}}\lesssim\left(J_{p,q}(\mu)\right)^{\frac{pq}{p-q}}.
$$
{{\em{Case 2}}}, $p=\frac{n}{\beta}:$ By the definition of
${\hbox{cap}}_{\dot{W}^{\beta,p}(\mathbb{R}^{n})}(O_{j}),$
 there is  $f_{j}\in \dot{W}^{\beta,p}(\mathbb{R}^{n})$ such that
 $$
 f_{j}\geq0, f_{j}(x)\geq1,\forall x\in O_{j} \ \hbox{and}\
\|f_{j}\|_{\dot{W}^{\beta,p}(\mathbb{R}^{n})}^{p}\leq2\
cap_{\dot{W}^{\beta,p}(\mathbb{R}^{n})}(O_{j})\leq 4\
c_{p}^{\beta}(\mu; 2^{j}).
$$
 Lemma \ref{le6} implies that  for each $j$ there is
$g_{j}(\cdot,\cdot)\in L^{p}(\mathbb{R}^{2n})$ such that
 $$
 f_{j}(x)=I_{2\beta}^{(2n)}\ast
 g_{j}(x,0)=\mathcal {R}\mathcal {E}f_{j}(x)
$$
and \begin{equation} \|I_{2\beta}^{(2n)}\ast
g_{j}\|_{\mathcal{L}^{p}_{2\beta}(\mathbb{R}^{2n})}=\|\mathcal{E}f_{j}\|_{\mathcal{L}^{p}_{2\beta}(\mathbb{R}^{2n})}
\leq\|f_{j}\|_{\dot{W}^{\beta,p}(\mathbb{R}^{n})}.\label{a6}
 \end{equation}
 Given integers $i, k$ with $i<k,$ define
  $$
  g_{i,k}=\sup_{i\leq j\leq
k}\left(\frac{2^{j}}{c_{p}^{\beta}(\mu;2^{j})}\right)^{\frac{1}{p-q}}g_{j}.
$$
Since $L^{p}(\mathbb{R}^{2n})$ is a lattice, we can conclude that
$g_{i,k}\in L^{p}(\mathbb{R}^{2n})$ and  $I_{2\beta}^{(2n)}\ast
g_{i,k}\in \dot{{\mathcal{L}}}_{2\beta}^{p}(\mathbb{R}^{2n}).$ Then
(\ref{a6}) and  Lemma \ref{le6} imply that
 \begin{eqnarray*}
&&  \|\mathcal{R}(I_{2\beta}^{(2n)}\ast
g_{i,k})\|^{p}_{\dot{W}^{\beta,p}(\mathbb{R}^{n})}\\
&\leq&
\sum\limits_{j=i}^{k}\left(\frac{2^{j}}{c_{p}^{\beta}(\mu;2^{j})}\right)^{\frac{p}{p-q}}
\|\mathcal{R}(I_{2\beta}^{(2n)}\ast
g_{j})\|_{\dot{W}^{\beta,p}(\mathbb{R}^{n})}^{p}
 \\ &\leq&
\sum\limits_{j=i}^{k}\left(\frac{2^{j}}{c_{p}^{\alpha}(\mu;2^{j})}\right)^{\frac{p}{p-q}}
\|I_{2\beta}^{(2n)}\ast
g_{j}\|_{\dot{\mathcal{L}}_{2\beta}^{p}(\mathbb{R}^{2n})}^{p}\\
&\lesssim&\sum\limits_{j=i}^{k}\left(\frac{2^{j}}{c_{p}^{\beta}(\mu;2^{j})}\right)^{\frac{p}{p-q}}
\|\mathcal{E}f_{j}\|_{\dot{\mathcal{L}}_{2\beta}^{p}(\mathbb{R}^{2n})}^{p}\\
&\lesssim&\sum\limits_{j=i}^{k}\left(\frac{2^{j}}{c_{p}^{\beta}(\mu;2^{j})}\right)^{\frac{p}{p-q}}
\|f_{j}\|_{\dot{W}^{\beta,p}(\mathbb{R}^{n})}^{p}\\
  &\lesssim&
\sum\limits_{j=i}^{k}\left(\frac{2^{j}}{c_{p}^{\beta}(\mu;2^{j})}\right)^{\frac{p}{p-q}}c_{p}^{\beta}(\mu;
2^{j}).\nonumber
  \end{eqnarray*}
  Note that for $i\leq j\leq k,$
$$
x\in O_{j}\Longrightarrow \mathcal{R}(I^{(2n)}_{2\beta}\ast
g_{i,k})(x)\geq \left(\frac{2^{j}}{c_{p}^{\beta}
(\mu;2^{j})}\right)^{\frac{1}{p-q}}.
$$
Then Lemma \ref{le4} (d) implies that
$$
(t,x)\in T(O_{j})\Longrightarrow S_{\alpha}(t^{2\alpha})|\mathcal
{R}(I^{(2n)}_{\beta}\ast g_{i,k})|(x)
\geq\left(\frac{2^{j}}{c_{p}^{\beta}(\mu;2^{j})}\right)^{\frac{1}{p-q}}\theta_{2}.
$$
This gives
$$
2^{j}<\mu\left(T(O_{j})\right)\leq\mu\left(E^{\alpha,\beta}_{\left(\frac{2^{j}}{c_{p}^{\beta}
(\mu;2^{j})}\right)^{\frac{1}{p-q}}\left(\frac{\theta_{2}}{2}\right)}(\mathcal
{R}(I^{(2n)}_{2\beta}\ast g_{i,k})(x))\right).
$$
 Hence
  \begin{eqnarray*}
&&\left(J_{p,q}(\mu)\|\mathcal{R}(I_{2\beta}^{(2n)}\ast
g_{i,k})\|_{\dot{W}^{\beta,p}(\mathbb{R}^{n})}\right)^{q}\nonumber\\
&\gtrsim&\!
\int_{\mathbb{R}_{+}^{1+n}}|S_{\alpha}(t^{2\alpha})\mathcal
{R}(I_{2\beta}\ast g_{i,k})(x)|^{q}\hbox{d}\mu(t,x)\\\nonumber
\!&\approx&\! \int_{0}^{\infty}\left(\inf\{\lambda:
\mu\left(E^{\alpha,\beta}_{\lambda}\mathcal {R}(I_{2\beta}\ast
g_{i,k})(x)\right)\leq s\}\right)^{q}\hbox{d}s\\\nonumber\!
&\gtrsim&\! \sum\limits_{j=i}^{k}\left(\inf\{\lambda:
\mu\left(E_{\lambda}^{\alpha,\beta}\mathcal {R}(I_{2\beta}\ast
g_{i,k})(x)\right)\leq 2^{j}\}\right)^{q}2^{j}\\\nonumber
\!&\gtrsim&\!\sum\limits_{j=i}^{k}
\left(\frac{2^{j}}{c_{p}^{\beta}(\mu;2^{j})}\right)^{\frac{q}{p-q}}2^{j}\\\nonumber
\!&\gtrsim&\!\!\frac{\sum_{j=i}^{k}\left(\frac{2^{j}}{c_{p}^{\beta}(\mu;2^{j})}\right)^{\frac{q}{p-q}}2^{j}}
{\left(\sum_{j=i}^{k}\left(\frac{2^{j}}{c_{p}^{\beta}(\mu;2^{j})}\right)^{\frac{q}{p-q}}c_{p}^{\beta}(\mu;
2^{j})\right)^{\frac{q}{p}}}\|\mathcal{R}(I_{2\beta}^{(2n)}\ast
g_{i,k})\|^{q}_{\dot{W}^{\beta,p}}\\\nonumber
\!&\approx&\!\!\left(\sum\limits_{j=i}^{k}\frac{2^{\frac{jp}{p-q}}}{\left(c_{p}^{\beta}(\mu;
2^{j})\right)^{\frac{q}{p-q}}}\right)^{\frac{p-q}{p}}\!\|\mathcal{R}(I_{2\beta}^{(2n)}\ast
g_{i,k})\|^{q}_{\dot{W}^{\beta,p}(\mathbb{R}^{n})}.\nonumber
  \end{eqnarray*}
 This tells us we obtain the same inequality as in the first case
 $$
 \sum\limits_{j=i}^{k}\frac{2^{\frac{jp}{p-q}}}{\left(c_{p}^{\beta}(\mu;
2^{j})\right)^{\frac{q}{p-q}}}\lesssim\left(J_{p,q}(\mu)\right)^{\frac{pq}{p-q}}.
$$
Note that the constant involved in the last inequality does not
depend on $i$ and $k.$ Letting $i\longrightarrow \infty$ and
$k\longrightarrow \infty,$ we have
$$
\int_{0}^{\infty}\left(\frac{t^{p/q}}{c_{p}^{\beta}(\mu;
t)}\right)^{\frac{q}{p-q}}\frac{dt}{t}
\lesssim\sum\limits_{-\infty}^{\infty}\frac{2^{\frac{jp}{p-q}}}{\left(c_{p}^{\beta}(\mu;
2^{j})\right)^{\frac{q}{p-q}}}\lesssim\left(J_{p,q}(\mu)\right)^{\frac{pq}{p-q}}.
$$
Therefore, (b) holds. $\Box$

 \vspace{0.1in}
\subsection{Proof of Theorem \ref{th3}} Let $p< q.$ The proof
consists two parts.

 {{{\it{Part 3.3.1}}}}: We prove $(a)\Longrightarrow(b)\Longrightarrow (c)\Longrightarrow(e)\Longrightarrow(a).$

$(a)\Longrightarrow(b)\Longrightarrow (c).$ Since
$\mu(E_{\lambda}(v_{0}))$ is nonincreasing in $\lambda,$
$$
q\mu(E^{\alpha,\beta}_{\lambda}(v_{0}))\lambda^{q-1}\leq\frac{d}{d\lambda}\left(\int^{\lambda}_{0}
\left(\mu(E_{s}^{\alpha,\beta}(v_{0}))\right)^{p/q}ds^{p}\right)^{q/p}.
$$
This gives,  for $s>0$
$$
(s^{q}\mu(E^{\alpha,\beta}_{s}(v_{0})))^{\frac{p}{q}}
\leq\left(q\int_{0}^{\infty}\mu(E^{\alpha,\beta}_{\lambda}(v_{0}))
\lambda^{q-1}\hbox{d}\lambda\right)^{\frac{p}{q}}
\leq\int_{0}^{\infty}\left(\mu(E^{\alpha,\beta}_{\lambda}(v_{0}))\right)^{\frac{p}{q}}\hbox{d}\lambda^{p},
$$
 and establishes the desired implications.\\  If  $(c)$ is true,   then
 $$
 K_{p,q}(\mu)\!=\!\sup\limits_{v_{0}\in\dot{W}^{\beta,p}(\mathbb{R}^{n}),\|v_{0}\|_{\dot{W}^{\beta,
p}}>0}
 \!\frac{\sup\limits_{\lambda>0}\lambda\left(\mu\left(\{(t,x)\in\mathbb{R}^{1+n}_{+}: |v(t^{2\alpha},x)|>\lambda\}
 \right)\right)^{\frac{1}{q}}}{\|v_{0}\|_{\dot{W}^{\beta,
p}}}<\infty.
$$
For a given $v_{0}\in \dot{W}^{\beta,p}(\mathbb{R}^{n})$ and a
bounded set $O\subseteq\hbox{Int}\left(\{x\in \mathbb{R}^{n}:
v_{0}(x)\geq 1\}\right),$ then Lemma \ref{le4} (d) implies
$$
\left(\mu\left(T(O)\right)\right)^{1/q}\lesssim
K_{p,q}(\mu)\|v_{0}\|_{\dot{W}^{\beta, p}}
$$
 and hence (e) follows from the definition of $cap_{\dot{W}^{\beta,p}}(O).$ To prove $(e)\Longrightarrow(a),$ we assume $(e).$ Then
$$
Q_{p,q}(\mu)=\sup\left\{\frac{\left(\mu\left(T(O)\right)\right)^{p/q}}
 {cap_{\dot{W}^{\beta,p}}(O)}:\ \hbox{bounded open }\ O\subseteq \mathbb{R}^{n}\right\}<\infty.
 $$
If $v_{0}\in\dot{W}^{\beta,p}(\mathbb{R}^{n})$ and $k=1,2,3,\cdots,$
then Lemmas \ref{le4} (a)-(c) and \ref{le1} (b) imply
 \begin{eqnarray*}
&&\int_{0}^{\infty}\left(\mu\left(E^{\alpha,\beta}_{\lambda}(v_{0})\cap
T (B(0,k))\right)\right)^{p/q}\hbox{d}\lambda^{p}\\\nonumber
 &\leq&\int_{0}^{\infty}\left(\mu\left(T(O^{\alpha,\beta}_{\lambda}(v_{0})\cap
B(0,k))\right)\right)^{p/q}\hbox{d}\lambda^{p}\\\nonumber
&\leq&\int_{0}^{\infty}\left(\mu_{k}\left(T(\{x\in\mathbb{R}^{n}:
\theta_{1}\mathcal{M}v_{0}(x)>\lambda\}\cap
B(0,k))\right)\right)^{p/q}\hbox{d}\lambda^{p}\\\nonumber
&\leq&\int_{0}^{\infty}\left(\mu\left(T(\{x\in\mathbb{R}^{n}:
\theta_{1}\mathcal{M}v_{0}(x)>\lambda\}\cap
B(0,k))\right)\right)^{p/q}\hbox{d}\lambda^{p}\\\nonumber
 & \lesssim&
Q_{p,q}(\mu)\int_{0}^{\infty}cap_{\dot{W}^{\beta,p}}\left(\{x\in\mathbb{R}^{n}:
\theta_{1}\mathcal{M}v_{0}(x)>\lambda\}\cap
B(0,k)\right)\hbox{d}\lambda^{p}\\\nonumber
 &\lesssim&
Q_{p,q}(\mu)\int_{0}^{\infty}cap_{\dot{W}^{\beta,p}}\left(\{x\in\mathbb{R}^{n}:
\theta_{1}\mathcal{M}v_{0}(x)>\lambda\}\right)\hbox{d}\lambda^{p}\\\nonumber
 &\lesssim& Q_{p,q}(\mu)\|v_{0}\|_{\dot{W}^{\beta, p}}^{p}.
\end{eqnarray*}

Letting $k\longrightarrow \infty$ in  the above inequality we have
$$
\int_{0}^{\infty}\left(\mu\left(E^{\alpha,\beta}_{\lambda}(v_{0})\right)\right)^{p/q}\hbox{d}\lambda^{p}
\lesssim Q_{p,q}(\mu)\|v_{0}\|_{\dot{W}^{\beta, p}}^{p}.
$$
 This derives (a).

{{{\it{Part 3.3.2}}}}: We verify
$(c)\Longrightarrow(d)\Longrightarrow(a).$

If (c) holds, then for any bounded open set
$O\subseteq\hbox{Int}\left(\{x\in \mathbb{R}^{n}: v_{0}(x)\geq
1\}\right),$ we have
 $$
 \mu\left(T(O)\right)^{1/q}\lesssim K_{p,q}(\mu)\|v_{0}\|_{\dot{W}^{\beta, p}}.
$$
Note that
 $$
 t^{p/q}\lesssim\left(K_{p,q}(\mu)\right)^{p}cap_{\dot{W}^{\beta,p}}(O)\ \hbox{wthenever}\ 0<t<\mu\left(T(O)\right).
$$
Hence
 $$
 t^{p/q}\lesssim \left(K_{p,q}(\mu)\right)^{p}c_{p}^{\beta}(\mu; t).
$$
 Therefore (d) holds.

If $(d)$ holds, then Lemmas \ref{le4} (b)-(c) and \ref{le1} (b)
imply that for each $k=1,2,3,\cdots,$
\begin{eqnarray*}
&&\vspace{0.1in}\int_{0}^{\infty}\left(\mu\left(E^{\alpha,\beta}_{\lambda}(v_{0})\cap
T(B(0,k))\right)\right)^{p/q}\hbox{d}\lambda^{p}\\
 &\leq&\int_{0}^{\infty}\left(\frac{\left(\mu\left(E^{\alpha,\beta}_{\lambda}(v_{0})\cap
T\left(B(0,k)\right)\right)\right)^{p/q}}{c_{p}^{\beta}\left(\mu;
\mu\left(E^{\alpha,\beta}_{\lambda}(v_{0})\cap
T\left(B(0,k)\right)\right)\right)}\right)cap_{\dot{W}^{\beta,p}}\left(O_{\lambda}(v_{0})\cap
 B(0,k)\right)\hbox{d}\lambda^{p}\\
&\lesssim& \left(\sup\limits_{t>0}\frac{t^{p/q}}{c_{p}^{\beta}(\mu;
 t)}\right)\int_{0}^{\infty}cap_{\dot{W}^{\beta,p}}\left(\{x\in \mathbb{R}^{n}: \theta_{1}\mathcal{M}
v_{0}(x)>\lambda\}\cap
 B(0,k)\right)\hbox{d}\lambda^{p}\\
 &\lesssim&\left(\sup\limits_{t>0}\frac{t^{p/q}}{c_{p}^{\beta}(\mu;
 t)}\right)\|v_{0}\|_{\dot{W}^{\beta,
p}}^{p} .
\end{eqnarray*}
Letting $k\longrightarrow \infty$ in  the previous inequality we
have
$$
\int_{0}^{\infty}\left(\mu\left(E^{\alpha,\beta}_{\lambda}(v_{0})\right)\right)^{p/q}\hbox{d}\lambda^{p}
\lesssim\left(\sup\limits_{t>0}\frac{t^{p/q}}{c_{p}^{\beta}(\mu;
 t)}\right)\|v_{0}\|_{\dot{W}^{\beta,
p}}^{p}.
$$
 This implies  that (a) holds. $\Box$

\vspace{0.1in} \subsection{Proof of Theorem \ref{th4}}

{{{\it{Part 3.4.1}}}}: We prove $(a)\Longleftrightarrow(b).$

It follows from Theorem \ref{th3} that it is enough to prove that
(b)
implies  (c) or (e) in Theorem \ref{th3}. We consider the following three cases.\\
{{{\em{Case 1}}}}, $1=p\leq q<\infty:$ If $(b)$ holds, then
$\|u\|_{1,q}<\infty.$ Suppose that $O\subseteq\mathbb{R}^{n}$ is a
bounded open set and is covered by a sequence of dyadic cubes
$\{I_{j}\}$ in $\mathbb{R}^{n}$ with
$\sum_{j}|I_{j}|^{\frac{n-\beta}{n}}<\infty.$ According to
Dafni-Xiao [8, Lemma 4.1] there exists another sequence of dyadic
cubes $\{J_{j}\}$ in $\mathbb{R}^{n}$ such that
$$
\hbox{Int}(J_{j})\cap\hbox{Int}(J_{k})=\emptyset\ \ \hbox{for}\ \
j\neq k,\ \ \bigcup_{j}J_{j}=\bigcup_{k}I_{k},
$$
$$
\sum_{j}|J_{j}|^{\frac{n-\beta}{n}}\leq\sum_{k}|I_{k}|^{\frac{n-\beta}{n}},\
T(O)\subseteq\bigcup_{j}T(\hbox{Int}(5\sqrt{n}J_{j})).
$$
Then
\begin{eqnarray*}
\mu(T(O))&\lesssim&\|\mu\|_{1,q}\sum\limits_{j}|5\sqrt{n}J_{j}|^{\frac{q(n-\beta)}{n}}\lesssim
\|\mu\|_{1,q}\left(\sum\limits_{j}|J_{j}|^{\frac{(n-\beta)}{n}}\right)^{q}\\
&\lesssim&\|\mu\|_{1,q}\left(\sum\limits_{j}|I_{j}|^{\frac{(n-\beta)}{n}}\right)^{q}.
\end{eqnarray*}
By  Xiao   \cite{J.Xiao 2}  (see also Adams \cite{D.R. Adams 2} or
\cite{D.R. Adams 3}) we have $cap_{1}^{\beta}(\cdot)\approx
H^{n-\beta}_{\infty}(\cdot),$ where the $H^{d}_{\infty}(\cdot)$ is
the $d-$ dimensional Hausdorff capacity. Thus, these along with the
definition of $H^{n-\beta}_{\infty}(O)$ imply
  $$
  \mu(T(O))\lesssim\|\mu\|_{1,q}\left(cap_{\dot{W}^{\beta,1}}(O)\right)^{q};
$$
 that is, the inequality (e) in Theorem \ref{th3} holds.

{{{\em{Case 2}}}}: $1< p<\min\{q, n/\beta\}:$ Let $v_{0}\in
\dot{W}^{\beta,p}(\mathbb{R}^{n})$ and
 $\mu_{\lambda}$ be the restriction of $\mu$ to $E^{\alpha,\beta}_{\lambda}(\dot{W}^{\beta,p}(\mathbb{R}^{n})).$
 If (b) holds, then
 $$
 \|\mu\|_{p,q}:=\sup\limits_{x\in\mathbb{R}^{n}, r>0}\frac{\mu\left(T(B(x,r))\right)}{r^{\frac{q(n-p\beta)}{p}}}<\infty.
 $$
It follows from Lemma \ref{le5} that
$$
|f(x)|\lesssim
\int_{\mathbb{R}^{n}}\frac{(-\triangle)^{\beta/2}f(y)}{|y-x|^{n-\beta}}\hbox{d}y,\
f\in \dot{W}^{\beta,p}(\mathbb{R}^{n}), x\in\mathbb{R}^{n}.
$$
 This inequality along with Lemma \ref{le2} and Fubini's theorem tell us
   \begin{eqnarray*}
\lambda\mu\left(E^{\alpha,\beta}_{\lambda}(v_{0})\right)\nonumber
&\lesssim&\int_{E^{\alpha,\beta}_{\lambda}(v_{0})}|S_{\alpha}(t^{2\alpha})(v_{0}(x))|\hbox{d}\mu(t,x)\\\nonumber
&\lesssim&\int_{E^{\alpha,\beta}_{\lambda}(v_{0})}\left|\int_{\mathbb{R}^{n}}K_{t^{2\alpha}}^{\alpha}(y)
(v_{0}(x-y))\hbox{d}y\right|\hbox{d}\mu(t,x)\\\nonumber
 &\lesssim&
\!\int_{\mathbb{R}^{1+n}_{+}}\!\left(\!\int_{\mathbb{R}^{n}}
\!\left(\!\int_{\mathbb{R}^{n}}\!\frac{K_{t^{2\alpha}}^{\alpha}(y)}{|z-(x-y)|^{n-\beta}}\hbox{d}y\!\right)
{\!|(-\triangle)^{\frac{\beta}{2}}v_{0}(z)|}\hbox{d}z\!\right)\hbox{d}\mu_{\lambda}(t,x)\\\nonumber
&\lesssim& \int_{\mathbb{R}^{1+n}_{+}}\left(\int_{\mathbb{R}^{n}}
\left((t^{2}+|z-x|^{2})^{\frac{\beta-n}{2}}\right){|(-\triangle)^{\beta/2}v_{0}(z)|}\hbox{d}z\right)\hbox{d}\mu_{\lambda}(t,x)\\\nonumber
&\lesssim&
\int_{\mathbb{R}^{n}}{|(-\triangle)^{\beta/2}v_{0}(z)|}\left(\int_{\mathbb{R}^{1+n}_{+}}
\left((t^{2}+|z-x|^{2})^{\frac{\beta-n}{2}}\right)\hbox{d}\mu_{\lambda}(t,x)\right)\hbox{d}z\\\nonumber
&\lesssim&
\int_{\mathbb{R}^{n}}{|(-\triangle)^{\beta/2}v_{0}(z)|}\left(\int_{0}^{\infty}
\mu_{\lambda}\left(T(B(z,r))r^{\beta-n-1}\right)\hbox{d}r\right)\hbox{d}z\\\nonumber
&\lesssim&(I_{1}(s)+I_{2}(s)),
  \end{eqnarray*}
 where
$$
I_{1}(s)=\int_{0}^{s}\left(\int_{\mathbb{R}^{n}}{|(-\triangle)^{\beta/2}v_{0}(z)|}
\mu_{\lambda}\left(T(B(z,r))\right)\hbox{d}z\right)r^{\beta-n-1}\hbox{d}r
$$
and
$$
I_{2}(s)=\int_{s}^{\infty}\left(\int_{\mathbb{R}^{n}}{|(-\triangle)^{\beta/2}v_{0}(z)|}
\mu_{\lambda}\left(T(B(z,r))\right)\hbox{d}z\right)r^{\beta-n-1}\hbox{d}r.
$$
By the definition of $\|\mu\|_{p,q},$ we have
$$
\mu_{\lambda}\left(T(B(z,r)\right)\leq\left(\mu_{\lambda}(T(B(z,r))\right)
^{1/p'}\|\mu\|_{p,q}^{1/p}r^{\frac{q(n-p\beta)}{p^{2}}}
$$
for $\frac{1}{p}+\frac{1}{p'}=1.$ So, using H$\ddot{o}$lder's
inequality and  the estimate
$$
\int_{\mathbb{R}^{n}}\mu_{\lambda}\left(T(B(x,r))\right)\hbox{d}x\lesssim
r^{n} \mu_{\lambda}\left(E^{\alpha,\beta}_{\lambda}(v_{0})\right),
$$
we obtain
  \begin{eqnarray*}
I_{1}(s)\nonumber&\lesssim&\!\int_{0}^{s}\left(\!\int_{\mathbb{R}^{n}}\!|(-\triangle)^{\beta/2}v_{0}(z)|
\!\left(\!\mu_{\lambda}(T(B(z,r)))\!\right)^{1/p'}\|\mu\|_{p,q}^{\frac{1}{p}}
r^{\frac{q(n-p\beta)}
{p^{2}}}\hbox{d}z\!\right)r^{\beta-n-1}\hbox{d}r\\\nonumber
 &\lesssim&\!\|v_{0}\|_{\dot{W}^{\beta,p}}\|\mu\|_{p,q}^{1/p}\int_{0}^{s}\left(\int_{\mathbb{R}^{n}}
 \mu_{\lambda}\left(T(B(z,r))\right)\hbox{d}z\right)^{1/p'}r^{\frac{q(n-p\beta)}{p^{2}}+\beta-n-1}\hbox{d}r\\\nonumber
 &\lesssim&\!\|v_{0}\|_{\dot{W}^{\beta,p}}\|\mu\|_{p,q}^{1/p}\int_{0}^{s}
 \left(r^{n}\mu\left(E^{\alpha,\beta}_{\lambda}(v_{0})\right)\right)^{1/p'}
 r^{\frac{q(n-p\beta)}{p^{2}}+\beta-n-1}\hbox{d}r\\\nonumber
&\lesssim&\!\|v_{0}\|_{\dot{W}^{\beta,p}}\|\mu\|_{p,q}^{1/p}
\left(\mu\left(E^{\alpha,\beta}_{\lambda}(v_{0})\right)\right)^{1/p'}s^{\frac{(q-p)(n-p\beta)}{p^{2}}}.\nonumber
  \end{eqnarray*}
 Similarly, we have
  \begin{eqnarray*}
  I_{2}(s)\!\!\nonumber&\lesssim&\!\!\int\limits_{s}^{\infty}\!\!\left(\int
\limits_{\mathbb{R}^{n}}\!{|(-\triangle)^{\frac{\beta}{2}}v_{0}(z)|^{p}}
\mu_{\lambda}\!\left(T(B(z,r))\right)\hbox{d}z\!\right)^{\frac{1}{p}}\!
\!\!\!\left(\int\limits_{\mathbb{R}^{n}}\!\!\mu_{\lambda}\!\left(T(B(z,r))\right)\hbox{d}z\!\right)
^{\frac{1}{p'}}\!\!\!\!\!\!r^{\beta-n-1}\!\hbox{d}r\\ \nonumber
 \!\!&\lesssim&\! \!\int_{s}^{\infty}\|v_{0}\|_{\dot{W}^{\beta,p}}
\left(\mu_{\lambda}\left(T(B(z,r))\right)\right)^{1/p}
\left(\int_{\mathbb{R}^{n}}\mu_{\lambda}\left(T(B(z,r))\right)\hbox{d}z\right)
^{1/p'}r^{\beta-n-1}\hbox{d}r\\ \nonumber \!
\!&\lesssim&\!\!\|v_{0}\|_{\dot{W}^{\beta,p}}
\left(\mu\left(E^{\alpha,\beta}_{\lambda}(v_{0})\right)\right)^{1/p}
\int_{s}^{\infty}r^{n/p'}\left(\mu(E^{\alpha,\beta}_{\lambda}(v_{0}))\right)^{1/p'}
r^{\beta-n-1}\hbox{d}r\\ \nonumber
 \!\!&\lesssim&\!\!\|v_{0}\|_{\dot{W}^{\beta,p}}
\left(\mu\left(E^{\alpha,\beta}_{\lambda}(v_{0})\right)\right)s^{\beta-n/p}.\nonumber
 \end{eqnarray*}
 Combing the above estimates on $I_{1}(s)$ and $I_{2}(s)$ together, we have
\begin{eqnarray*}
\lambda\mu(E^{\alpha,\beta}_{\lambda}(v_{0}))&\lesssim&\|v_{0}\|_{\dot{W}^{\beta,p}}
\mu_{\lambda}\left(E^{\alpha,\beta}_{\lambda}(v_{0})\right)\\
&\times&\left(s^{\beta-n/p}\!+\!
\left(\|\mu\|_{p,q}\left(\mu\left(E^{\alpha,\beta}_{\lambda}(v_{0})\right)\right)^{-1}\right)^{1/p}
s^{\frac{(q-p)(n-p\beta)}{p^{2}}}\right).
\end{eqnarray*}
 Taking
 $$
 s=\left(\|\mu\|_{p,q}^{-1}\left(\mu\left(E^{\alpha,\beta}_{\lambda}
 (v_{0})\right)\right)\right)^{\frac{p}{q(n-p\beta)}}
 $$
 in the above inequality, we have
$$
\lambda\left(\mu(E^{\alpha,\beta}_{\lambda}(v_{0}))\right)^{1/q}
\lesssim\|\mu\|_{p,q}^{1/q}\|v_{0}\|_{\dot{W}^{\beta, p}}.
$$
This implies  the condition (c) of Theorem \ref{th3}.\\

 {{{\it{Part 3.4.2}}}}: We find a nonnegative Radon measure  to show that if $1<p=q<n/\beta$ then (b) does not imply
 (a) in general.

 In fact, suppose  $K\subseteq\mathbb{R}^{n}$ is a compact set with  the $(n-p)-$dimensional Hausdorff measure
 $H^{(n-p\beta)}(K)>0,$ then by  Maz'ya \cite[p. 358, Proposition 3]{V.G. Maz'ya 1} we have $cap_{\dot{W}^{\beta,p}}(K)=0,$
 on the other hand by Adams-Hedberg  \cite[p. 132, Proposition 5.1.5 \&  p. 136, Theorem 5.1.12]{D.R. Adams Hedberg} we can
  find a nonnegative Radon  measure $\nu$ on   $\mathbb{R}^{n}$ such that
$$
\sup\limits_{x\in\mathbb{R}^{n},
r>o}\frac{\nu\left(B(x,r)\right)}{r^{n-p\beta}}<\infty\ \hbox{and}\
0<H^{n-p\beta}_{\infty}(K)\lesssim \nu(K).
$$
Define $\mu(t,x)=\delta_{1}(t)\otimes \nu(x).$ Then (b) hold for
this nonnegative  Radon  measure on $\mathbb{R}^{1+n}_{+}.$ However,
(a) is  not  true, otherwise, we would have $0<\nu(K)\lesssim
cap_{\dot{W}^{\beta,p}}(K)=0.$ Contradiction. $\Box$

\vspace{0.1in} \subsection{Proof of Theorem \ref{th5}}
Suppose $0<q<
1.$ Since the proof of $(a)\Longrightarrow(b)\Longrightarrow(c)$ is
similar
 to that of $(b)\Longrightarrow(c)\Longrightarrow(e)$  of Theorem \ref{th3}, we only need  to verify $(c)\Longrightarrow(d).$ Let $(c)$
  be true. Then Lemma \ref{le4} (a)$-$(c) imply
 \begin{eqnarray*}
 \mu\left(E^{\alpha,\beta}_{\lambda}(v_{0})\right)
 &\leq& \left(\mu\left(T(O^{\alpha,\beta}_{\lambda}(v_{0}))\right)\right)\\
&\leq& \left(\mu\left(T(\{x\in\mathbb{R}^{n}:
\theta_{1}\mathcal{M}v_{0}(x)>\lambda\})\right)\right)\\
 & \lesssim&
\left(cap_{\dot{W}^{\beta,1}}\{x\in\mathbb{R}^{n}:
\theta_{1}\mathcal{M}v_{0}(x)>\lambda\}\right)^{q}.
\end{eqnarray*}
 This and Lemma \ref{le1} (b) imply that (d) holds. $\Box$

\vspace{0.1in} \subsection{Proof of Theorem \ref{th6}} From the
proof of Theorems \ref{th3} \& \ref{th4} for  $1\leq p<n/\beta$ and
$q>p,$ we have
 $$
 \begin{array}{llll}
\||\mu|\|_{p, q}\!&=&\!\sup\limits_{x\in\mathbb{R}^{n}, r>0}
 \frac{\left(\mu\left(T(B(x,r))\right)\right)^{\frac{p}{q}}}{cap_{\dot{W}^{\beta,p}}\left(B(x,r)\right)}<\infty\\
&\Rightarrow&\left(\int\limits_{\mathbb{R}^{1+n}_{+}}|v(t^{2\alpha},x)|^{q}\hbox{d}\mu(t,x)\right)^{\frac{1}{q}}
\lesssim \||\mu|\|_{p,q}\|v_{0}\|_{\dot{W}^{\beta,p}}, \ \forall
v_{0}\in \dot{W}^{\beta,p}(\mathbb{R}^{n}).
  \end{array}
 $$
Given $(t_{0}, x_{0})\in \mathbb{R}^{1+n}_{+}.$ Let
$q=\frac{np}{n-p\beta}$ and $\mu(t,x)=\delta_{(t_{0},x_{0})}$ be the
Dirac measure at $(t_{0},x_{0}).$ It suffices to prove
$\||\delta_{(t_{0},x_{0})}|\|_{p,q}\leq t_{0}^{p\beta-n}.$
 In fact,  if $(t_{0},x_{0})$ is not in  $T(B(x,r)),$ then  $\delta_{(t_{0},x_{0})}(T(B(x,r)))=0.$ If $(t_{0},x_{0})\in T(B(x,r)),$ then
$B(x_{0},t_{0})\subseteq B(x,r)$ and $r^{n}\geq t_{0}^{n}.$ This
give the estimate
$$
\delta_{(t_{0},x_{0})}(T(B(x,r)))\leq
\frac{r^{n}}{t_{0}^{n}}=t_{0}^{-n}r^{\frac{(n-p\beta)q}{p}}.
$$
The above estimate  and $cap_{\dot{W}^{\beta,p}}(B(x,r))\approx
r^{n-p\beta}$ verify
$$
\frac{(\delta_{(t_{0},x_{0})}(T(B(x,r))))^{p/q}}{cap_{\dot{W}^{\beta,p}}(B(x,r))}\leq
t_{0}^{-\frac{np}{q}}.
$$
 Therefore, $\||\delta_{(t_{0},x_{0})}|\|_{p,q}\leq t_{0}^{p\beta-n}.$ $\Box$

\vspace{0.08in} \subsection{Proof of Theorem \ref{th9}} Assume that
$\mu$ is a nonnegative Radon  measure such  that
$$
\sup\limits_{x\in \mathbb{R}^{n},
r>0}\frac{\left(\mu\left(T(B(x,r))\right)\right)^{p/q}}
{cap_{\dot{W}^{1/2,p}}(B(x,r))}<\infty
 $$
 for  $1\leq p<2n$  and $\frac{4pn+4p}{2n-p}\leq q<\infty.$ According to the definition of  $1/2-$parabolic rectangle
  (see Nishio-Yamada \cite{M. Nishio})
 $$
 Q^{1/2}(s,y)=\{(s,y)\in\mathbb{R}^{1+n}_{+}: |x_{j}-y_{j}|<s/2, 1\leq j\leq n, s\leq t\leq 2s\}
 $$
  with center $(s,y),$  we have
  $$
  Q^{1/2}(s,y)=[s, 2s]\times B(y, \sqrt{n}s/2).
  $$
 The definition of $T(B(y, r))$ implies  that there is  a  dimensional constant $c(n)$ such that
$$
 Q^{1/2}(s,y)\subseteq  T(B(y, c(n)s)),
$$
 for each $(s,y)\in \mathbb{R}^{1+n}_{+},$ so
$$
\mu(Q^{1/2}(s,y))\leq \mu(T(B(y, c(n)s)))\lesssim s^{q(n-p/2)/p}.
$$
If $p_{1}=\frac{q(2n-p)}{2p(n+1)}-1,$   then for each $(s,y)\in
\mathbb{R}^{1+n}_{+},$
$$
\mu(Q^{1/2}(s,y))\lesssim s^{(n+1)(1+p_{1})}.
$$
Note that $p_{1}\geq1$ since $q\geq \frac{4p(n+1)}{2n-p}$ and
$p<2n.$ It follows from Nishio-Yamada \cite[p. 91 Theorem 2]{M.
Nishio} that $\nu$ is a (0,1)-type  Carleson  measure on
$b_{1/2}^{q}$ $(q\geq 1)$ if and only if $\nu(Q^{1/2}(s,y))\lesssim
s^{(n+1)(1+q)},$  for each $(s,y)\in \mathbb{R}^{1+n}_{+}.$ Thus
$\mu$ is a (0,1)-type Carleson
 measure on $b_{1/2}^{p_{1}}.$ $\Box$

   \vspace{0.1in}
\noindent
 {\bf{Acknowledgement.}}

 This work is a part of my doctoral thesis. I want to thank my supervisor Professor Jie Xiao
  for suggesting the problem
 with kind  encouragement.
 I am  grateful to the referee for several helpful suggestions on the manuscript,
especially for  pointing out the results proved by
Cascan-Ortega-Verbitsky in  \cite{Cascan ortega verbitsky}.


\begin{thebibliography}{1}

\bibitem{D.R. Andams 76}
D.R. Adams,
 \textit{On the existence of capacitary strong type
estimates in $\mathbb{R}^{n},$} {Ark. Mat.} 14 (1976) 125-140.

\bibitem{D.R. Adams 4}
 D.R. Admas, \textit{Lectures on $L^{p}-$potentia Theory}, Univ. of Umea, 1981.

\bibitem{D.R. Adams 2}
D.R. Adams,
 \textit{The classification problem for the capacities
associated with the  besov and Triebel$-$Lizorkin spaces},
 {Banach Center Publ}.
 22 (1989) 9-24.

 \bibitem{D.R. Adams 3}
 D.R. Adams,
\textit{Choquet integrals in potential theorey}, {Publ. Mat.} 42
(1998) 3-66.



\bibitem{D.R. Adams Hedberg}
 D.R. Adams, L.I. Hedberg,
\textit{Function Spaces and Potential Theory}, {A Series of
Comprehensive Studies in Mathematics, Springer, Berlin}, 1996.

\bibitem{D. R. Adams J. Xiao}
D.R. Adams, J. Xiao,
 \textit{ Strong type estimates for homogeneous
Besov capacities}, {Math. Ann.} 325 (2003) 695-709.





\bibitem{L. Carleson}
L. Carleson,
 \textit{Interpolation by bounded analytic functions
and the Corona problem}, {Ann. of Math.} 76 (1962) 547-559.

\bibitem{C. CASCANTE J.M. ORTEGA AND I.E. VERBITSKY 1}
C. Cascante, J.M. Ortega, I.E. Verbitsky,
 \textit{Trace inequalities of
Sobolev type in the upper triangle case,} {Proc. London Math. Soc.,}
80 (2000) 391¨C414.


\bibitem{Cascan ortega verbitsky}
C. Cascante, J.M., Ortega, I.E., Verbitsky,
 \textit{Wolff's inequality for radially nonincreasing kernels and applications to trace
inequalties}, {Potential Anal.} 16 (2002) 347-372.

\bibitem{Costea Maz'ya}
S. Costea, V.G. Maz'ya,
 \textit{Conductor inequalities and criteria for Sobolev-Lorentz two-weight inequalities}, {preprint,
arXiv:0804.3051}.



\bibitem{Dafni karadzhov Xiao}
G. Dafni, G. E. Karadzhov, J. Xiao,
 \textit{Classes of Carleson-type measures generated by capacities}, {Math. Z.} 258 (2008) 827-844.


\bibitem{Dafni Xiao 2}
 G. Dafni, J. Xiao,
  \textit{Some new tent spaces and duality theorems for fractional Carleson measures and $Q_\alpha(\Bbb R^n)$,} {J.
Funct. Anal.}, 208 (2004) 377-422.


 \bibitem{L. Grafakos}
L. Grafakos,
 \textit{Classical and Modern Fourier Analysis}, {Praeson Education, Inc.}, 2004.

\bibitem{W. Hastings}
W. Hastings,
 \textit{A Carleosn measure theorem for Bergman spaces}, {Proc. Amer. Soc.}, 52 (1975) 237-241.

\bibitem{R. Johnson}
 R. Johnson,
 \textit{Application of Carleson measures to partial differential equations and Fourier multiplier problems, in:  Proceedings
 of the Conference on harmonic Analysis},
 {Cortona, lecture Notes in Mathematics, Vol. 992, Springer, Berlin},  16-72.

\bibitem{Karadzhov Xiao}
G. E. Karadzhov, J. Xiao,
 \textit{Carleson type theorems for certain convolution operators}, {Integr. Equ. Oper. Theory} 55 (2006)
429-438.


\bibitem{V.G. Maz'ya 62}
V.G. Maz'ya,
  \textit{The negative spectrum of the $n-$dimensional Schr$\ddot{o}$dinger operator.} {Dokl. Akad. Nauk SSSR,} 144 (1962)
721-722 (Russian). {English translation: Soviet Math. Dokl.} 3
(1962) 808-810.

\bibitem{V.G. Maz'ya 64}
V.G. Maz'ya,
  \textit{On the theory of the  multidimensional Schr$\ddot{o}$dinger operator.} {Izv.  Akad. Nauk SSSR,} 28 (1964)
1145-1172 (Russian).


\bibitem{V.G. Maz'ya 2}
V.G. Maz'ya,
 \textit{On certain integral inequalities for functions of many varibles}, {J. Sov. Math.} 1 (1973) 205-234.

\bibitem{Maz'ya}
V.G. Maz'ya,
 \textit{On capacity strong type estimates for fractional norms}. {Zup. Nauch. Sem. Leningrad otel. math. Inst.
Steklov} (LOMI) 70 (1977) 161-168, (Russian).

\bibitem{V.G. Maz'ya 3}
V.G. Maz'ya,
 \textit{Estimates for capacities and traces of potentials},
  {Internat. J. Math. $\&$  Math. Sci.} 7 (1984) 41-63.


\bibitem{V.G. Maz'ya 1}
V.G. Maz'ya,
 \textit{Sobolev Spaces}, {Springer, Berlin, New York,} 1985.

\bibitem{V.G. Maz'ya 4}
V.G. Maz'ya,
 \textit{Lectures on isoperimetric and isocapacitary inequalities in the theory of Sobolev spaces},  {Heat kernels and
analysis on manifolds, graphs, and metric spaces} (Paris, 2002),
 307--340, {Contemp. Math., 338, Amer. Math. Soc., Providence, RI}, 2003.

\bibitem{Maz'ya 9}
V. Maz'ya,
 \textit{Conductor and capacitary inequalities for functions on topological spaces and their applications to Sobolev
 type imbeddings,} {Journal of Functional Analysis,} 224 (2005) 408-430.


\bibitem{Maz'ya 5}
V.G. Maz'ya,
 \textit{Integral and isocapacitary inequalities,} {preprint, arXiv:0809.2511}


\bibitem{V.G. Maz'ya Havin}
 V.G. Maz'ya, V.P. Havin,
\textit{Nonlinear potential theory}, {Usp. Mat. Nauk}, 27 (1972)
 67-138 (Russian). {English translation: Russian Math. Surveys} 27 (1972) 71-148.

\bibitem{Mazya Netrusov}
V.G. Maz'ya, Y. Netrusov,
 \textit{Some counterexamples for the theory of
 Sobolev spaces on bad domians,} {potential Analysis,} 4 (1995) 47-65.


\bibitem{Maz'ya Preobra}
V.G. Maz'ya, S.P. Preobra$\breve{z}$enski$\breve{{\i}},$
 \textit{On estimates of $(p,l)-$capacity and traces of potentials.} {Wissenschaftliche Informationen. Technische hochschule,
Karl-Marx-Stadt, Sektion Mathematik,} 28 (1981) 1-38 (Russian).

\bibitem{Mazya shaposhnikova}

V.G. Maz'ya, T.O. Shaposhnikova, \textit{The Theorey of Multipliers
in spaces of Differentiable Functions}, {Monographs and Studies in
Mathematics}, 23, Pitman, Boston-Londan, 1985.

\bibitem{Maz'ya Verbitsky 1}
V. G. Maz'ya, I.E. Verbitsky,
 \textit{Capacitary estimates for fractional integrals, with applications to partial differential equations and
Sobolev multipliers,} {Arkiv f$\ddot{o}$r Matem.} 33 (1995) 81-115.


\bibitem{C. Miao}
C. Miao, B. Yuan, B. Zhang,
 \textit{Well-posedness of the Cauchy problem for the fractional power dissipative equations},  {Nonlinear Anal. TMA}
  68 (2008) 461-484.

\bibitem{M. Nishio K. Shimomura N. Suzuki}
M. Nishio, K. Shimomura, N. Suzuki,
 \textit{$\alpha-$parabolic Bergman spaces}, {Osaka J. Math.} 42 (2005)  133-162.


\bibitem{M. Nishio}
 M. Nishio, M. Yamada,
 \textit{ Carlson type measures on parabolic Bergman spaces},  {J. Math. Soc. Japan} 58 (2006) 83-96.





\bibitem{E.M. Stein}
E.M. Stein,
 \textit{ Harmonic Analysis: Real-varible Methods, Orthogonality, and
 Oscillatory Integrals}, {Princeton University
Press, Princeton, New Jersey}, 1993.

\bibitem{H. Triebel}
H. Triebel,
 \textit{Theory of Function Spaces II}, {Birkh$\ddot{a}$user,
Basel}, 1992.

\bibitem{Verbitsky}
I. E. Verbitsky,
 \textit{Nonlinear potentials and trace
inequalities,} {The Maz¡¯ya Anniversary Collection, Eds. J.
Rossmann, P. Tak¡äac, and G. Wildenhain, Operator Theory: Advances
and Applications} 110 (1999), 323-343.


 \bibitem{Z. Wu}
 Z. Wu,
 \textit{Strong type estimate and Carleson measures for Lipschitz spaces},
 {Proc. Amer. Math. Soc.} 127 (1999) 3243-3249.


\bibitem{J. Xiao}
J. Xiao, \textit{Carleson embeddings for Sobolev spaces via heat
equation}, {J. Differ. Equations} 224 (2006) 277-295.

\bibitem{J.Xiao 2}
J. Xiao, \textit{Homogeneous endpoint Besov space embeddings by
 hausdorff capacity and heat equation},{ Adv. in Math.} 207 (2006) 828-846.

\bibitem{J. Xiao 3}
J. Xiao,
 \textit{The $Q_{p}$ Carleson measure problem}, {Adv. in Math.} 217 (2008) 2075-2088.

\end{thebibliography}
\end{document}